\documentclass[12pt]{amsart}

\usepackage{amsmath}
\usepackage{amssymb}
\usepackage{fancybox}
\usepackage{eucal}
\usepackage{amsthm}
\usepackage{amscd}
\usepackage{wasysym}
\usepackage{url}
\usepackage{MnSymbol} %symbols for infinitary formulas

\usepackage[all,cmtip]{xy}

\usepackage{enumitem}% http://ctan.org/pkg/enumitem

\usepackage{hyperref}

\usepackage[dvips]{graphicx}
\usepackage{latexsym}
\usepackage{mathrsfs} 

\newtheorem{thm}{Theorem}[section]
\newtheorem{prop}[thm]{Proposition}

\newtheorem{lemma}[thm]{Lemma}

\newtheorem{cor}[thm]{Corollary}

\newtheorem{rmk}[thm]{Remark}

\newtheorem{definitiontemp}[thm]{Definition}
\newenvironment{defn}{\begin{definitiontemp}
\normalfont}{\end{definitiontemp}}

\newcommand{\Fields}{\operatorname{\textup{\textbf{Fields}}}}

\newcommand{\DCFields}{\operatorname{\textup{\textbf{DCFields}}}_0}

\newcommand{\DCFC}{\operatorname{\textup{\textbf{DCFC}}}}

\newcommand{\DCFT}{\operatorname{\textup{\textbf{DCFC}}}}

\newcommand{\Graphs}{\operatorname{\textup{\textbf{Graphs}}}}

\newcommand{\DD}{\mathscr{D}}
\newcommand{\FF}{\mathscr{F}}
\newcommand{\GG}{\mathscr{G}}
\newcommand{\HH}{\mathscr{H}}
\newcommand{\DDtilde}{\widetilde{\mathscr{D}}}
\newcommand{\GGtilde}{\widetilde{\mathscr{G}}}

\newcommand{\CC}{\mathbf{C}}

\newcommand{\TKbar}{\overline{T_K}}
\newcommand{\TKFbar}{\overline{T_{K_F}}}
\newcommand{\TKHbar}{\overline{T_{K_H}}}
\newcommand{\TKH}{T_{K_H}}
\newcommand{\TKRbar}{\overline{T_{K_R}}}
\newcommand{\TQbar}{\overline{T_{\Q}}}
\newcommand{\TQA}{T_{\Q\la A\ra}}
\newcommand{\TQAbar}{\overline{T_{\Q\la A\ra}}}

\newcommand\aut[1]{\text{Aut}(#1)}
\newcommand\Aut[1]{\text{Aut}(#1)}

\newenvironment{pf}{\begin{trivlist}\item[\hskip\labelsep
{\it Proof.}]}{\end{trivlist}}

\newcommand{\comment}[1]{}

\newcommand{\Iso}[1]{\text{Iso}(#1)}

\newcommand{\ACF}[1]{\ensuremath{\textbf{ACF}_{#1}}}
\newcommand{\DCF}[1]{\ensuremath{\textbf{DCF}_{#1}}}

\newcommand{\set}[2]{\ensuremath{ \{ #1 : #2 \} }}

\newcommand{\spec}[1]{\ensuremath{ \text{Spec} (#1 )}}

\renewcommand{\deg}[1]{\ensuremath{\text{deg}(#1)}}

\newcommand{\Z}{\mathbb{Z}}

\newcommand{\Q}{\mathbb{Q}}
\newcommand{\QA}{\Q\la A\ra}

\newcommand{\C}{\mathcal{C}}
\newcommand{\A}{\mathcal{A}}
\newcommand{\B}{\mathcal{B}}

\newcommand{\Kbar}{\overline{K}}

\newcommand{\QAbar}{K_0}

\newcommand{\Atilde}{\widetilde{A}}

\newcommand{\Gtilde}{\widetilde{G}}
\newcommand{\Htilde}{\widetilde{H}}

\newcommand{\atilde}{\tilde{a}}

\newcommand{\ptilde}{\tilde{p}}
\newcommand{\qtilde}{\tilde{q}}

\newcommand{\utilde}{\tilde{u}}
\newcommand{\vtilde}{\tilde{v}}

\newcommand{\xtilde}{\tilde{x}}
\newcommand{\ytilde}{\tilde{y}}

\newcommand{\Ktilde}{\widetilde{K}}

\newcommand{\K}{\mathcal{K}}

\newcommand{\M}{\mathfrak{M}}

\renewcommand{\SS}{\mathcal{S}}

\newcommand{\avec}{\vec{a}}

\newcommand{\Uvec}{\vec{U}}
\newcommand{\Vvec}{\vec{V}}
\newcommand{\uvec}{\vec{u}}
\newcommand{\vvec}{\vec{v}}

\newcommand{\xvec}{\vec{x}}
\newcommand{\Xvec}{\vec{X}}
\newcommand{\yvec}{\vec{y}}
\newcommand{\Yvec}{\vec{Y}}

\newcommand{\del}{\delta}

\newcommand{\KY}{K\{ Y\}}

\newcommand{\la}{\langle}
\newcommand{\ra}{\rangle}
\def\embeds{\hookrightarrow}

\newcommand{\at}{\char'100}

\newcommand{\proves}{\vdash}

\newcommand{\dom}[1]{\text{dom}(#1)}

\newcommand{\bfd}{\boldsymbol{d}}
\newcommand{\Qhat}{\widehat{\mathbb{Q}}}

\newcommand{\khat}{\widehat{k}}
\newcommand{\Khat}{\widehat{K}}
\newcommand{\KFhat}{\widehat{K_F}}
\newcommand{\KHhat}{\widehat{K_H}}

\def\bfz{\boldsymbol{0}}
\def\s01{\ensuremath{\Sigma^0_1}}
\def\d02{\ensuremath{\Delta^0_2}}
\def\phi{\varphi}

\def\res{\!\!\upharpoonright\!}

\begin{document}

\title[Differentially Closed Fields and Universality on a Cone]{Differentially Closed Fields\\and Universality on a Cone}

\author[R. Miller]{Russell Miller}
\address{Dept.\ of Mathematics, Queens College, 65-30 Kissena Blvd.,
Queens, NY 11367 \& Ph.D.\ Programs in Mathematics \& Computer Science,
Graduate Center, City University of New York, 365 FIfth Avenue, New York, NY 10016, USA.}
\email{Russell.Miller@qc.cuny.edu}
\urladdr{\href{http://qcpages.qc.cuny.edu/~rmiller/}{http://qcpages.qc.cuny.edu/$\sim$rmiller/}}

\thanks{
The research for this article was partially supported by Grant \#581896 from
the Simons Foundation, and by several grants from
The City University of New York PSC-CUNY Research Award Program.
It benefitted from useful conversations with David Marker, Matthew
Harrison-Trainor, and Tom Scanlon while the author participated in a program supported
by the National Science Foundation under Grant \# DMS-1928930 and hosted by the
Simons Laufer Mathematical Sciences Institute in Berkeley, California, during the summer of 2022.
%This article expands upon a talk given (remotely) by the author in the
%Algebra \& Logic Seminar in Novosibirsk on 23 September 2021.
}

\maketitle

\begin{abstract}
The class of all countable differentially closed differential fields $K$ of
characteristic $0$ was shown by Marker and the author in \cite{MM18} to be ``one jump away''
from universal for spectra of structures:  for every nontrivial countable
structure $\M$, there is some $K$ whose spectrum is the preimage
under jump of the spectrum of $\M$, and conversely, for every $K$,
there is such an $\M$.  We show that
the missing jump can be accounted for by adding to the signature
of differential fields a predicate describing a certain algebraic transcendence property.
The ensuing universality results for differentially closed fields in the
new signature include not only spectra of structures, but also many
properties related to computable categoricity.  However, these latter
universality results hold only on the cone above a specific $\Sigma^0_1$
oracle set, whose decidability status remains unknown.  Moreover,
differentially closed fields simply fail flat-out to be universal for automorphism
groups, even non-effectively.

We also include a small erratum (Remark \ref{rmk:correction}) to \cite{MM18}.
\end{abstract}

\section{Introduction}
\label{sec:intro}

The theory \DCF0 of ordinary differentially closed fields of characteristic $0$ has fundamental
similarities to the theory \ACF0 of algebraically closed fields of that characteristic.
Both theories are complete and decidable, with quantifier elimination.
Model theorists, however, have identified several respects in which \DCF0
is more complex than \ACF0.  For example, \DCF0 has Morley rank
$\omega+1$ (see \cite[Cor.\ 5.14]{M06}), whereas for \ACF0, the Morley rank is $1$:
every definable subset of an algebraically closed field must be either finite or cofinite,
whereas an algebraic differential equation can have infinitely many solutions.
%(Below we describe some further differences.)

In \ACF0, of course, countable models are characterized by their transcendence
degrees, which must be cardinals $\leq\omega$, giving only countably
many such models.  Hrushovski and Sokolovi\'c \cite{hs} established another
significant difference between these theories
by showing that \DCF0 has continuum-many countable models.
To accomplish this, they constructed, for each countable graph $G$,
a model $\DD(G)$ of \DCF0 such that $\DD(G)\cong \DD(\Gtilde)$ if and only if
$G\cong\Gtilde$.  (\cite{P06} describes their result nicely.)
In \cite{MM18}, Marker and the author examined
the extent to which their construction of $\DD(G)$ from $G$ is effective.
It turns out to be Turing-computable:  from the atomic diagram $D(G)$,
one can produce the atomic diagram $D(\DD(G))$ uniformly and effectively.
Indeed, more is true:  $D(\DD(G))$ is uniformly computable even from
an approximation of $D(G)$.  In the opposite direction, however, this fails:
it is not generally possible to compute the atomic diagram of a copy of $G$
from that of an arbitrary copy of $\DD(G)$.  The main result of \cite{MM18}
is expressed as follows.  (Recall the definition:  for a countable structure $\A$,
$\spec{\A}=\set{\deg{D(\B)}}{\B\cong\A~\&~\dom{\B}=\omega}$.)
\begin{thm}[Theorem 5.5 of \cite{MM18}]
\label{thm:MM}
A set $\SS$ of Turing degrees is the spectrum of a model of \DCF0 if and only if
there exists a countable, automorphically nontrivial graph $G$ such that
$$ \SS = \set{\bfd}{\bfd'\in\spec{G}}.$$
\end{thm}
In particular, for a graph $G$, $\spec{\DD(G)}$ is the preimage of
$\spec{G}$ under the jump operation $\bfd\mapsto \bfd'$.

The class $\Graphs$ of all symmetric, irreflexive, automorphically nontrivial graphs
on the domain $\omega$ is known to be \emph{universal} for all standard computable-structure-theoretic
properties.  This result appears mainly in \cite{HKSS02}, though various aspects
of it had been proven earlier.  It has also turned out to hold for certain properties
that were unknown when \cite{HKSS02} first appeared.
The term \emph{complete} is often used interchangeably
with ``universal'' to name this property.  We choose ``universal'' here because we will use
``complete'' in several other senses:  to refer to maximal consistent theories,
maximal consistent types, formulas $\phi(x)$ that generate complete types,
graphs with no missing edges, and subsets of $\omega$ of maximal complexity
within their own level in the arithmetical or the Borel hierarchy.

Universality includes the fact that
every spectrum of an automorphically nontrivial structure is also
the spectrum of some graph in $\Graphs$.
It includes a substantial further list of properties as well, which will
be detailed in Section \ref{sec:consequences}.
Our purpose in this article is to make precise the intuition from Theorem \ref{thm:MM}
that the class $\DCFields$ of countable models of \DCF0 is ``one jump away''
from being universal.  We will show that, when the signature of differential fields
is augmented by a unary relation symbol $C$ that holds of precisely those
elements algebraic over a particular differential subfield $\QA$ of the
differential closure $\Qhat$ of the rational numbers, the resulting class
$\DCFC$ of countable models of \DCF0 in this signature is universal
for many computable-model-theoretic properties,
in the same sense that $\Graphs$ is, except that the universality
holds \emph{on a cone}.  Below we will elaborate on this concept.
Of course, $C(x)$ is a $\Sigma_1$ property of $x$, defined by the computable
$L_{\omega_1\omega}$ formula that states that $x$ is a solution
to some nontrivial algebraic equation over $\QA$.  (Moreover,
$A$ itself is definable without parameters within $\Qhat$.) 
So the one jump by which $\DCFields$ was ``off'' is accounted for
by the difficulty of deciding algebraicity over $\QA$ in models of \DCF0.
On the other hand, we will note in Section \ref{sec:notafunctor} that
these universality results cannot be extended to automorphism groups:
in that context, graphs simply have more complexity than models of \DCF0.

For clarity:  $C$ defines field-theoretic algebraicity over $\QA$,
not differential algebraicity.  In Section \ref{sec:binarydependence}
we discuss the alternative of replacing $C$
with a more naturally defined symbol, namely
the binary relation of algebraic dependence over $\Q$ itself,
with equivalent results.

In \cite{MPSS18}, Poonen, Schoutens, Shlapentokh,
and the author used the concept of a \emph{computable functor}
to demonstrate that the class $\Fields_p$ of countable
fields (of any fixed characteristic $p\geq 0$) is universal in the sense of \cite{HKSS02}.
Of course, functors arise mainly in category theory, and indeed those authors
considered $\Fields_p$ as a category, exactly as we will do here for the classes
named above.
\begin{defn}
\label{defn:cats}
$\Fields_p$, $\DCFields$, $\DCFT$, and $\Graphs$ are the categories
in which the objects are all those structures with domain $\omega$ in the given class
(respectively: fields of characteristic $p\geq 0$, models of the theory \DCF0,
models of \DCF0 in the signature with the algebraicity predicate $C$,
and symmetric irreflexive automorphically nontrivial graphs)
and the morphisms from an object $M$ to another object $N$
are precisely those isomorphisms mapping $M$ onto $N$.
(The elements of $\Graphs$ may be described as the symmetric
infinite coinfinite subsets of $(\omega^2-\set{(x,x)}{x\in\omega})$.)
\end{defn}
Since the algebraicity predicate $C$ is $L_{\omega_1\omega}$-definable,
the categories $\DCFC$ and $\DCFields$ have exactly the same objects
and exactly the same morphisms.  However, the same object will have
different atomic diagrams in the two categories (one with $C$ and one without),
and for purposes of computability, these are distinct categories.
Specifically, a presentation $K\in\DCFields$ may have computable atomic
diagram $D(K)$ there, yet as an object in $\DCFC$, the same $K$ may no longer
be computable, since algebraicity over $\QA$ in $K$ may be undecidable.
It should also be noted that, unlike $\DCFields$, $\DCFC$ is not the set
of all countable models of a complete theory, as it is not possible for a theory
to ensure that $C$ holds of precisely those elements algebraic over $\QA$.
(There is an infinite axiom scheme stating that $C$ holds of \emph{all}
such algebraic elements, but the converse requires an infinitary formula.)

For our purposes here, it is necessary both to relativize and to generalize
the notion of a computable functor, which was first defined in \cite{MPSS18}.
First we relativize.
\begin{defn}
\label{defn:functor}
Let $S\subseteq\omega$.  An \emph{$S$-computable functor} $\FF$ from any of the
categories $\CC$ above to any other one $\mathbf{D}$ is a functor in the usual sense of category theory,
with the additional property that there exist Turing functionals $\Phi$ and $\Phi_*$ such that,
for all structures $M$ and $N$ in $\CC$ and for all isomorphisms $f:M\to N$ in $\CC$,
\begin{itemize}
\item
$\Phi^{S\oplus D(M)}$ computes the atomic diagram of $\FF(M)$; and
\item
$\Phi_*^{S\oplus D(M) \oplus f\oplus D(N)} = \FF(f)$, computing the isomorphism $\FF(f)$
from the domain $\omega$ of $\FF(M)$ onto the domain $\omega$ of $\FF(N)$.
\end{itemize}
If the fixed oracle set $S$ is decidable, then the functor is a \emph{computable functor}
in the sense of \cite[Defn.\ 3.1]{MPSS18}.
\end{defn}

Originally it was our hope to apply this notion to $\Graphs$ and $\DCFT$,
using the construction from \cite{MM18}, in the same way that \cite{MPSS18}
did for $\Graphs$ and $\Fields$.  This is not possible, for reasons that will be discussed below.
Instead, we will use the weaker notion of a \emph{reduction of categories},
which deletes the requirement usually called ``functoriality.''

\begin{defn}
\label{defn:reduction}
A \emph{reduction of categories} $\FF$ from any of the
categories $\CC$ above to any other one $\mathbf{D}$ consists of two maps
(both named $\FF$, as for functors), one mapping each object $M$ of $\CC$
to an object $\FF(M)$ of $\mathbf D$, and the other mapping each morphism
$f:M\to N$ of $\CC$ to a morphism $\FF(f):\FF(M)\to\FF(N)$ in $\mathbf D$.
The map $\FF$ on morphisms is \emph{not} required to preserve composition,
nor to map identity morphisms to identity morphisms.

Fix $S\subseteq\omega$.  A reduction of categories $\FF$ is \emph{$S$-computable}
if it has the additional property that there exist Turing functionals $\Phi$ and $\Phi_*$ such that,
for all structures $M$ and $N$ in $\CC$ and for all isomorphisms $f:M\to N$ in $\CC$,
\begin{itemize}
\item
$\Phi^{S\oplus D(M)}$ computes the atomic diagram of $\FF(M)$; and
\item
$\Phi_*^{S\oplus D(M) \oplus f\oplus D(N)} = \FF(f)$, computing the isomorphism $\FF(f)$
from the domain $\omega$ of $\FF(M)$ onto the domain $\omega$ of $\FF(N)$.
\end{itemize}
If the oracle set $S$ is decidable, then $\FF$ is a \emph{computable reduction of categories}.
\end{defn}

We will imitate the example of \cite{MPSS18}, using an oracle
to compute a reduction $\DD$ from $\Graphs$ into $\DCFC$
and then producing a computable reduction $\GG$ from the image of $\DD$
back into $\Graphs$, such that the compositions in both directions are
the identity functors on those categories (up to a \emph{computable natural isomorphism}, as defined below).
The first twist is that $\DD$ is not a functor, although $\GG$ will be one.
This will eliminate certain universality properties from consideration:
notably, our construction does not establish the universality
of $\DCFC$ for automorphism groups.

The second twist is that, while $\GG$ will be computable, $\DD$ will only be
computable relative to an oracle set $T_\infty$, which contains certain formulas that
are complete in a particular theory $\DCF0\cup D(K_\infty)$.  In Subsection \ref{subsec:oracle}
we discuss the particular model $K_\infty$ of \DCF0 and this oracle set defined using it.
$T_\infty$ is a specific subset of $\omega$, defined in differential algebra, whose
decidability status is presently an open question:  it is known to be a
$\Sigma^0_1$ set, but, as far as is known, its Turing degree could be
any c.e.\ degree, from $\bfz$ up to $\bfz'$.  If $T_\infty$ should turn out to be decidable,
then both reductions of categories would be computable, and the universality
results for $\DCFC$ in Section \ref{sec:consequences} would hold in full.
If not, then $\DCFC$ may be the first example of a natural class of structures
that is universal relative to an oracle, but not universal in its own right.
In the latter case, we would say that $\DCFC$ is \emph{universal on
the cone above $T_\infty$}.  This will be explained in Section \ref{sec:consequences}.
Regardless of which of these cases holds, however, the results will reinforce
the conclusion of \cite{MM18}
that differentially closed fields have far greater computable-structure-theoretic
complexity than algebraically closed fields.

\section{Background on differential fields}
\label{sec:diffalg}
\subsection{Formulas generating types}
\label{subsec:types}

The differential closure $\Khat$ of a differential field $K$ is defined to be the prime model
of the theory $\textbf{DCF}\cup D(K)$, just as the algebraic closure of a field $F$
is the prime model of $\textbf{ACF}\cup D(F)$.  In both cases the characteristic is
determined by the atomic diagram, $D(K)$ or $D(F)$; here we will only consider fields
of characteristic $0$.  As a prime model, the differential closure realizes precisely
the principal types over this theory $\textbf{DCF}\cup D(K)$ (equivalently, $\DCF0\cup D(K)$),
and so each element $a\in\Khat$ is characterized (up to automorphisms over $K$)
by some formula $\alpha(x)$ such that $\Khat\models\alpha(a)$ and such that
$\alpha(x)$ generates an entire (complete) type over $\textbf{DCF}\cup D(K)$.
A formula $\alpha(x)$ with this property is said to be a \emph{complete formula}
for the theory \DCF0.  Since \DCF0 has quantifier elimination, one may assume
this generating formula $\alpha(x)$ to be quantifier-free.

The same analysis holds in $\textbf{ACF}\cup D(F)$, and there the generating formula
may always be taken to be of the form $p(x)=0$, where $p$ ranges over
the irreducible polynomials in $F[X]$.  For certain fields $F$, irreducibility
of polynomials in $F[X]$ can fail to be decidable from $D(F)$, but for all finitely
generated fields, it is decidable.  In differential fields, the most natural analogues
of these generating formulas are the \emph{constrained pairs} $(p,q)$ from
the ring $K\{ Y\}=K[Y,Y',Y'',\ldots]$ of differential polynomials over $K$.
\begin{defn}
\label{defn:constrainedpair}
Let $p,q\in K\{Y\}$.  The pair $(p,q)$ is \emph{constrained} if:
\begin{itemize}
\item
$p$ is monic and irreducible and has strictly larger order than $q$; and
\item
for every $h\in K\{Y\}$ of lower order than $p$, and for every $y,z\in\Khat$
satisfying $p(y)=0\neq q(y)$ and $p(z)=0\neq q(z)$,
$$ h(y)=0 \iff h(z)=0.$$
\end{itemize}
We write $\TKbar = \set{(p,q)\in (K\{Y\})^2}{(p,q)\text{~ is constrained}}$ for the
set of all constrained pairs over $K$.  The complementary set $T_K$ of unconstrained pairs
is computably enumerable uniformly in $D(K)$.
\end{defn}
In order to ensure that each $a\in K$ realizes a constrained pair,
we adopt the convention that nonzero constant polynomials have order $-1$:
thus, even though $Y-a$ has order $0$, the pair $(Y-a,1)$ lies in $\TKbar$.
%and gives the complete formula $Y-a=0\neq 1$.

Using quantifier elimination, we may conclude that for each constrained
pair $(p,q)$, the formula $p(Y)=0\neq q(Y)$
generates a type.  Conversely, it is known that every principal type over
$\DCF0 \cup D(K)$ is generated by such a formula.  So these formulas naturally
play the same role as the irreducible polynomial equations over a field $F$.
(Indeed, the formulas $p(Y)=0\neq q(Y)$ are essential to the standard axiomatization
of \DCF0, which was given by Blum in \cite{B68}.)

However, as remarked earlier, \DCF0 is in some ways more complex than \ACF0,
and at present it remains open whether $T_{\Q}$ is decidable.  Here $\Q$ is viewed
as a \emph{constant} differential field, in which every element has derivative $0$:
thus $\Q$ is the prime model of the theory $\textbf{DF}_0$ of differential fields
(not necessarily differentially closed!) of characteristic $0$.

For clarity we remark here that all theorems stated here about principal types
and complete formulas apply in the standard signature of differential fields.
When the relation $C$ is adjoined to that signature, many things change.
Adjoining $C$ does not change any isomorphism relations, as $C$ is
$L_{\omega_1\omega}$-definable, but since it is not definable by any finitary
formula, one cannot readily transfer theorems about types from $\DCF0$
into the larger signature.  However, in a $K\in\DCFT$, each element still realizes
a principal type in the signature of $\DCF0$, and every such type is still realized;
only the formulas involving $C$ cause problems.

In \cite{H74}, Harrington showed that every computable differential field $K$ of characteristic
$0$ has a computable differential closure.
In doing so, he showed that there is a uniform computable enumeration of all the principal
types over the theory $\DCF0\cup D(K)$.  However, as he remarked himself, his proof
does not show how to identify, within this enumeration, a generator of each type.
His proof relativizes readily to yield the following result.
\begin{cor}[cf.\ Corollary 3(i) in \cite{H74}]
\label{cor:Harrington}
There are Turing functionals $\Theta$ and $\Gamma$ such that, whenever $D(K)$ is the atomic diagram
of a differential field $K$ of characteristic $0$ (coded as a subset of $\omega$, by a G\"odel coding),
$\Theta^{D(K)}$ computes the atomic diagram of a differentially closed field $\Khat$
and $g_K=\Gamma^{D(K)}$ computes a differential field embedding of $K$ into $\Khat$
such that $\Khat$ is the differential closure of the image $g_K(K)$ of this embedding.
Thus $\Khat$ realizes all principal types of $\DCF0\cup D(g_K(K))$ but realizes no other types.
\qed\end{cor}

\subsection{Useful Fact}
\label{subsec:tools}

\begin{prop}
\label{prop:alg}
Let $K$ be a differential field, with differential closure $\Khat$.
Then the algebraically closed subfield
$$\Kbar= \set{y\in\Khat}{y\text{~is algebraic over~}K}$$
of $\Khat$ is in fact a differential subfield, i.e., closed under the differential operator $\del$,
and there is a Turing functional $\Phi$ that, given the atomic diagram
of $K$ as an oracle, produces the atomic diagram of (a copy of) $\Kbar$
as a differential field, computing derivatives as well as $+$ and $\cdot$,
and also produces a differential field embedding $f:K\embeds\Kbar$
such that $\Kbar$ is algebraic over the image $f(K)$.
\end{prop}
\begin{pf}
Closure of $\Khat$ under the differential operator $\del$ is immediate:  if $y\in\Khat$ has
$h(y)=\sum a_ny^n=0$ with all $a_n\in K$, then
$$ 0=\del(h(y))=\sum \del(a_n)y^n + \del(y)\cdot\sum na_ny^{n-1} ,$$
giving us a formula for $\del(y)$ as an (algebraic) rational function of $y$ over $K$.

Rabin's Theorem gives an effective presentation $\Kbar$ of the algebraic closure
of $K$, along with the necessary embedding $f:K\embeds\Kbar$.
Having the formula above (as opposed to a mere polynomial in $K[T]$ with root $\del(y)$)
allows us to compute $\del$ and thus extend the atomic diagram of $\Kbar$
to include atomic facts about $\del$.  Thus we can avoid the use of Harrington's
Theorem, although it would also suffice.  Notice that the computation of $\del(y)$
above does not require the algebraic polynomial $h\in K[Y]$ to be irreducible --
which is important, because we are not assuming that $K$ has a splitting algorithm.
\qed\end{pf}

Clearly $\Kbar$ is computably enumerable within the differential closure $\Khat$
given by Harrington's Theorem.  One might hope for $\Kbar$ to be decidable there.
This holds in some simple cases, notably when $K$ itself is a constant field (as then
$y\in\Kbar\iff\del(y)=0$).  In general, however, this would require decidability of the set
$T_K$ of complete formulas for $K$:  with a $T_K$-oracle, one could find the constrained
pair $(p,q)$ satisfied by a given $y\in\Khat$, and conclude that $y\in\Kbar$
just if $p$ has order $0$.  This situation will be important in the constructions below.

\section{Construction of the Reductions $\DD$ and $\GG$}
\label{sec:functors}

As promised, we now construct the two reductions of categories necessary for the results
in Section \ref{sec:consequences}.  The first one, $\DD$, maps the category
$\Graphs$ into the category $\DCFT$ of differentially closed fields of characteristic $0$,
in the signature with the algebraicity predicate $C$ defining algebraicity
over $\QA$ (as described below).  This $\DD$
will require an oracle set; its (near-)inverse $\GG$ will be a computable functor.

The procedure $\DD$ on a graph $G$ in $\Graphs$ is best understood by imagining
the edges of $G$ to be enumerated, rather than decided.  To formalize this,
we describe quickly a second pair of computable functors $\HH$ and $\FF$,
to go back and forth between these concepts.  The output computed
by $\HH$ on input $D(G)$ will be an \emph{enumeration} of the edges in the graph
$\HH(G)$, i.e., a subset of $\omega^3$ whose image under the projection 
$(x,y,s)\mapsto (x,y)$ is the set of all pairs of adjacent nodes in $\HH(G)$.

$\HH$ accepts as input the atomic diagram $D(G)$
of a graph $G$ on the domain $\omega$.  For each node $x\in G$, the new graph
$\HH(G)$ has a node $c_x$, identifiable because a loop of three other nodes
is also added to the graph, with one of the three adjacent to $c_x$.  For every $x<y$,
we also add a node $d_{xy}$, adjacent to both $c_x$ and $c_y$, and we attach to $d_{xy}$
either a loop of length $5$ (if $D(G)$ says that $x$ and $y$ are adjacent)
or a loop of length $7$ (if $D(G)$ says that they are not.)  This defines the graph
$\HH(G)$, but the output of $\HH$ on input $G$ is an enumeration of the
edges in this new graph.

The inverse functor $\FF$ accepts as input any enumeration of the edges
in a graph $H$ isomorphic to any graph in the image of $\HH$.  Using this enumeration,
$\FF$ eventually identifies each node $c$ adjacent to a loop of length $3$ in $H$
and creates a node $x_c$ in $\FF(H)$.  Give any $x_c$ and $x_{c'}$, it then
finds the unique node $d$ in $H$ adjacent to both $c$ and $c'$, waits until it sees
a loop of length either $5$ or $7$ appear in $H$ adjacent to this $d$,
and outputs accordingly (in the atomic diagram of $\FF(H)$) whether
$x_c$ and $x_{c'}$ have an edge between them or not.

These processes are clearly both effective.  
Moreover, given an isomorphism $g:G\to\Gtilde$, $\FF$ can compute the
obvious isomorphism $\FF(g)$ from $\FF(G)$ onto $\FF(\Gtilde)$, and $\HH$
can do likewise with an isomorphism from $H$ onto $\Htilde$.
$\FF\circ\HH$ and $\HH\circ\FF$ are not actually the
identity functors, but they are effectively isomorphic to them.
Since $\HH$ and $\FF$ essentially just formalize the notion of Marker $\exists$-extensions,
we do not feel compelled to give any further details here.  Below,
Proposition \ref{prop:DcircG} will add a few specifics.

Now we can turn to our version of the procedure used in \cite{MM18}, which is the map
on objects used by the reduction $\DD$ from $\Graphs$ to $\DCFT$.
Given a graph $G$, it first applies $\HH$ to $D(G)$ to produce an enumeration
of $H=\HH(G)$.  Next, it takes a fixed computable differential ground field which we
will call $\QAbar$, isomorphic to the algebraic closure of the differential subfield
$\QA$ generated within (a computable presentation of) $\Qhat$ by the subset
defined by the \emph{Rosenlicht equation}:
$$A=\set{y\in\Qhat}{y\neq 0~\&~y\neq 1~\&~y'=y^3-y^2},$$
which is known to be a strongly minimal set of indiscernible elements within $\Qhat$.
(The computable presentation $\Qhat$ is given by the theorem of Harrington
from \cite{H74}.  $\QAbar$ is a computably enumerable subfield of $\Qhat$ --
indeed a differential subfield, by Proposition \ref{prop:alg} --
hence computably presentable in its own right.)
Writing $A=\{ a_0<a_1<\ldots\}$, our procedure treats each $a_n$ as the
representative of the node $n\in H$.  If the enumeration of edges in $H$
ever indicates that there is an edge between the nodes $m$ and $n$
(with $m<n$), our procedure adjoins to the differential field
a new pair of elements $(u_{mn}, v_{mn})$.  Each coordinate individually
is transcendental over $\QA$, but the two together satisfy
the elliptic-curve equation for $(m,n)$:  we set
$$ v_{mn}^2 = u_{mn}(u_{mn}-1)(u_{mn}-a_m-a_n).$$
It is well known that the solutions to the equation $y^2=x(x-1)(x-a_m-a_n)$
form an abelian group.  $\QAbar$ already contains many pairs satisfying this equation,
-- indeed, for each $j>0$,
$K_0$ contains exactly $j^2$ elements of order $j$ in the group -- but
all solutions in $\QAbar$ have coordinates algebraic over $\QA$.  In contrast,
the newly added pair will have infinite order and each of its coordinates
is transcendental over $\QA$.  This pair with transcendental coordinates
codes the existence of the edge between $m$ and $n$ into the differential field we are building.

(To be clear:  each of $u_{mn}$
and $v_{mn}$ has minimal differential polynomial of order $1$ over $\QA$.
These elements are differentially algebraic but algebraically trancendental
over $\QA$.  Together they are algebraically dependent over $\QA$,
since $ v_{mn}^2 = u_{mn}(u_{mn}-1)(u_{mn}-a_m-a_n)$.  If the group element
$(u_{mn},v_{mn})$ had had finite order, that would have yielded another algebraic
equation satisfied by the element, saying that its order was $j$, in which case
each coordinate would have become algebraic over $\QA$.  Instead, each coordinate
in the new point realizes a non-principal type in \DCF0 over $\QA$, and therefore
this type is not realized in $\Qhat$.)

$\DD$ thus constructs a differential field extension of $\QA$:
$$K_H = \QAbar \la u_{mn},v_{mn}~:~m<n~\&~(m,n)\in H\ra.$$
\begin{lemma}
\label{lemma:splittingalg}
$\DD$ can compute a presentation of $K_H$, with a splitting algorithm,
uniformly in the graph enumeration $H=\HH(G)$.
Moreover, each of the subfields $\Q$, $\QA$, and $K_0$ will be decidable within $K_H$,
uniformly in the enumeration of $H$, and will have its own splitting algorithm.
\end{lemma}
\begin{pf}
Computing the atomic diagram of $K_H$ only requires starting with a computable
copy of $K_0$ and then, each time a new edge $(m,n)$ appears in $H$,
adjoining the new $u_{mn}$ and $v_{mn}$ to this differential field.
All these $u_{mn}$ have the same minimal differential polynomial,
except that each uses its own $a_m$ and $a_n$ in that polynomial.
(The fact that $u_{mn}$ has this same minimal differential polynomial
even over the field generated by the preceding $u_{m'n'}$ and $v_{m'n'}$
follows from orthogonality of the types of the pairs $(u_{mn},v_{mn})$,
as in the original Hrushovski-Sokolovi\'c construction.)
Notice that the minimal differential polynomial of $u_{mn}$ over $K_0$
(or equivalently, over the extension of $K_0$ by other $u_{m'n'}$'s)
gives a transcendence basis for $K_0\la u_{mn}\ra$, namely $\{ u_{mn},u_{mn}',
\ldots,u_{mn}^{(r-1)}\}$, where $r$ is the order of that polynomial.  Additionally,
the polynomial itself serves as the minimal polynomial of $u_{mn}^{(r)}$
over $K_0(u_{mn},\ldots,u_{mn}^{(r-1)})$, and $\{u_{mn},\ldots,u_{mn}^{(r)}\}$
generates $K_0\la u_{mn}\ra$ as a field over $K_0$.
Finally, the elliptic curve equation gives the minimal polynomial
of $v_{mn}$ over $K_0\la u_{mn}\ra$, again without interference
from the other solutions already adjoined.

Since the procedure above is uniform and the minimal polynomials are all known,
a splitting algorithm for $K_H$ can be computed uniformly from one for $K_0$.
As $K_0$ is the algebraic closure of $\QA$, this only requires knowing a transcendence
basis for the field $\QA$ over $\Q$, and $A$ itself is such a basis (as every element
of $A$ has order $1$ over $\Q$) and is decidable, being defined by the Rosenlicht
equation.  As a field, $\QA$ is generated by $A\cup\set{a'}{a\in A}$, and the minimal
polynomial of each $a'$ over the purely transcendental extension $\Q(A)$
is the Rosenlicht polynomial, so the existence of a splitting algorithm for $\Q$
(which is the original fact proven by Kronecker) yields a splitting algorithm
for $\QA$ as well.

With these splitting algorithms and bases, we get the decidability promised in the Lemma.
For any $y\in K_H$, find its minimal (algebraic) polynomial over the subfield
$K_0(u_{mn}~:~(m,n)\in H)$, expressed as $p\in K_0[u_{m_0n_0},\ldots,u_{m_kn_k},X]$ 
for some finite collection of basis elements $u_{mn}$.  This will determine whether
$x\in K_0$ or not.  If $x\in K_0$, we can go further and find its minimal polynomial
over $\QA$ (which will determine whether it lies in $\QA$) and also over $\Q(A)$
(which will determine whether it lies in $\Q$).
\qed\end{pf}

This extension $K_H$ is not differentially closed, but $\DD$ applies
Corollary \ref{cor:Harrington} to build its differential closure $\KHhat=\DD(G)$
as it goes along, making $\DD(G)$ indeed a model of \DCF0.
It was shown by Hrushovski and Sokolovi\'c in \cite{hs} that,
for those $(m,n)$ such that no edge between $m$ and $n$ ever appears in $H$,
$\KHhat$ will contain no pair $(x,y)$ of points transcendental over $\QA$
satisfying the elliptic curve equation for $a_m$ and $a_n$.
%Equivalently, $\KHhat$ contains no non-torsion points in the relevant abelian variety.
(This result also appears in \cite[Section 5]{M00} and \cite[Section 4]{P06}.)
Moreover, neither the adjoinment of transcendental points representing edges in $H$
nor the subsequent extension to the differential closure $\KHhat$ adds any
new elements satisfying the definition of the set $A$, so
the set thus defined within $\KHhat$ is the same set $A$.
Thus the coding of $H$ into $\DD(G)$ is successful -- and so, in turn,
is the coding of $G$.

\begin{rmk}[An erratum.]
\label{rmk:correction}
In \cite{MM18}, a more involved construction succeeds in computing
a copy of $\DD(G)$ just from an approximation of the atomic diagram
of $G$.  The authors there, including the present author, exploited the fact
that elements that appear to be transcendental over $\QA$ can later be made
algebraic (and turned into torsion points of the abelian variety in question)
if the approximation so dictates.  In doing so, they misstated one
aspect of the situation:  they claimed that all solutions $(x,y)$
to the elliptic curve equation that have $x$ and $y$ algebraic
over $\QA$ are torsion points for this abelian variety.  This holds
for certain other elliptic curves (e.g., for $y^2=x(x-1)(x-a)$ with $a\in A$,
because each $a\in A$ has order $1$ over $\Q$).  However,
$(a_m+a_n)$ has order $2$ over $\Q$, and it is an open question
whether the claim holds for the elliptic curves $y^2=x(x-1)(x-a_m-a_n)$.
Nevertheless, regardless of the status of the claim in that case,
the results in \cite{MM18} do hold:  when the authors speak
of nontorsion points over $y^2=x(x-1)(x-a_m-a_n)$, one need only assume
that they mean nontorsion points with coordinates transcendental over
$\QA$.  (As there are no torsion points with transcendental coordinates,
it would be equivalent simply to say ``solutions $(x,y)$ with $x$ and
$y$ each transcendental over $\QA$.'')  The properties of being nontorsion
and having transcendental coordinates are both definable by computable
$\Pi_1$ formulas of $L_{\omega_1\omega}$, so the arguments in \cite{MM18}
adapt easily to this slightly different characterization and the results there still hold.
\end{rmk}

Here we do not attempt the more involved construction of \cite{MM18}:
we started with $\QAbar$ as the ground field, with all torsion
points already present there (possibly along with some nontorsion points with
coordinates algebraic over $\QA$, as the existence of such points in $\Qhat$ remains open).
Instead, we need to decide the predicate $C$ defining algebraicity over $\QA$,
since $\DD$ is to be a reduction to $\DCFC$, not just to $\DCFields$.
With $D(G)$ (and thus the enumeration of $H$) available to it,
$\DD$ can do this using the oracle set $T$ described below.
Corollary \ref{cor:oracle} below shows that $T_{K_H}$ is computable
uniformly from $T\oplus D(G)$.  Knowing $T_{K_H}$, we can now
take any $y\in\KHhat=\DD(G)$ and find the unique $(p,q)\in\overline{T_{K_H}}$
such that $p(y)=0\neq q(y)$.  If $p(Y)$ has positive order, then
certainly $y$ is not algebraic over $\QA$.  If $p(Y)$ is an algebraic
polynomial (i.e., of order $0$), then $y$ is algebraic over $\QA$
just if all the coefficients of $p$ lie in $\QA$:  indeed $p\in K_H[Y]$
and the field extension $K_H/\QAbar$ is purely transcendental.
Thus, modulo the proof of Corollary \ref{cor:oracle},
$\DD$ does produce the atomic diagram of $\DD(G)$
in the signature for $\DCFC$, as required.  This is the procedure
from Subsection \ref{subsec:tools}.
Notice that it would not have been possible to decide $C$ without access
to the enumeration of $H$:  an approximation to $G$ (and thus a
$\Sigma_2$ presentation of the edges in $H$) would fall one jump short.

To establish that $\DD$ is a reduction of categories, we must also compute
$\DD(g)$ whenever $g:G_0\to G_1$ is an isomorphism of graphs.
This will be done below.  First, having explained the construction
of $\DD(G)$, we describe the inverse functor $\GG$, which accepts as input
the atomic diagram of an element $K$ of $\DCFT$ that is isomorphic
to a differential field in the image of $\DD$.  Knowing the basic structure
of $K$, therefore, it begins by going through the elements of $K$
and identifying the elements $a_0<a_1<\cdots$ used to code the graph,
which are precisely those elements $y$ satisfying the formula
$y\neq 0~\&~y\neq 1~\&~y'=y^3-y^2$ defining $A$.  %(Notice that the predicate
%$T$ allows us to decide membership in $\Qhat$.  However,
%$\Qhat$ would be enumerable from $D(K)$ in any case, so this is not yet crucial.)
For each $a_n$, $\GG$ creates a corresponding node $x_n$ in the graph
$H=\GG(K)$.  Next, for each $m<n$, it watches for a pair $(u,v)$ to appear
in $K$ satisfying $v^2=u(u-1)(u-a_m-a_n)~\&~\neg C(u)%~\&~\neg C(v)
$.  Here the predicate
$C$ is essential, because many pairs of elements algebraic over $\QA$ also satisfy
this equation.  If such a pair $(u,v)$ of transcendentals ever appears,
then $\GG$ enumerates an edge into $H$ between $x_m$ and $x_n$.
Finally, $\GG$ uses the functor $\FF$ above to convert this enumeration
of the edges in $H$ into a computation of the atomic diagram of the graph
$\FF(H)$, which is the final output $\GG(K)$.

Once again this computation differs by a jump from that in \cite{MM18}.
Here we have $C$, which allows $\GG$ to recognize immediately
any pair that codes the existence of an edge in $H$.  In \cite{MM18},
it was necessary to guess whether such a pair was transcendental or not
(equivalently, whether it had infinite or finite order in the abelian variety),
and so that construction gave the edge relation in $H$ as a $\Sigma_2$ set,
and thus gave only an approximation to the atomic diagram of $\DD(K)$.

The computation of isomorphisms by $\GG$ is straightforward, and requires
no oracle.  Given as input an isomorphism $f:K\to\Ktilde$ (along with
the atomic diagrams of $K$ and $\Ktilde$), $\GG$ readily recognizes
the elements $a_0<a_1<\cdots$ of $A$ in $K$ and $\atilde_0<\atilde_1<\cdots$
of $\Atilde$ in $\Ktilde$.  Each of these sets is defined (within $K$ and $\Ktilde$)
by the same formula, so $f$ must map $A$ bijectively onto $\Atilde$.
If $f(a_m)=\atilde_n$, then $\DD(g)$ must map $x_m$ to $\xtilde_n$,
in the graphs $H$ and $\Htilde$ enumerated by the first part of the process.
This is clearly an isomorphism from $H$ onto $\Htilde$, since $K$
will contain a transcendental solution to $y^2=x(x-1)(x-a_m-a_n)$
if and only if $\Khat$ contains such a solution to $y^2=x(x-1)(x-f(a_m)-f(a_n))$.
Finally, $\GG$ uses $\FF$ to construct the isomorphism from the graph
$G=\FF(H)=\GG(K)$ onto $\Gtilde=\FF(\Htilde)=\GG(\Ktilde)$,
and this is the output $\GG(f)$.  It is clear that this result is functorial:
$\GG(f_0\circ f_1)=\GG(f_0)\circ\GG(f_1)$, and applying $\GG$
to the identity map on any $K$ will yield the identity map on $\GG(K)$.

The remaining step is the most interesting:  showing how $\DD$
accepts an isomorphism $g:G\to\Gtilde$ of graphs and computes
an isomorphism $f=\DD(g):\DD(G)\to\DD(\Gtilde)$.  Intuitively it is
clear what the procedure should do.  First it should compute the isomorphism
$h=\HH(g):H\to\Htilde$, where $H=\HH(G)$ and $\Htilde=\HH(\Gtilde)$.
(This $\HH$ is functorial.)
Next, since nodes in $H$ are represented by elements of $A$ within
the subfield $\Qhat$ of $K=\DD(G)$, it should use $h$ to guide its
choice of a bijection from $A$ onto $\Atilde\subseteq\Qhat\subseteq\Ktilde=\DD(\Gtilde)$.
In particular, if $h(m)=n$, then $f(a_m)$ should equal $\atilde_n$.
As $h:H\to\Htilde$ is a bijection, so will $f\res A$ be.  Moreover,
as $A$ is a set of indiscernibles in $\Qhat$, every permutation of $A$
extends to an automorphism of $\Qhat$ and hence to an automorphism
of the algebraic closure $K_0$ of $\QA$ within $\Qhat$.
So $f$ should apply this
automorphism to $K_0$ (details below!), mapping the subfield $K_0$ of $K$
to the subfield $K_0$ of $\Ktilde$.  ($K$ and $\Ktilde$ were each built
around the same fixed computable presentation $K_0$, so the automorphism
of $K_0$ may be regarded as an isomorphism from the subfield $K_0$ of $K$
onto the subfield $K_0$ of $\Ktilde$).  This defines $f\res K_0$.
Next, whenever new elements $(u,v)$ were added to $K$ in
the computation of $K=\DD(G)$, they were added as a point with transcendental
coordinates in the abelian variety defined by $y^2=x(x-1)(x-a_m-a_n)$ for some
$m$ and $n$, and this was done because an edge appeared between
$m$ and $n$ in the enumeration of edges in $H$.  Since $g$ is an isomorphism,
so is $h=\HH(g)$, so there must also be an edge between $h(m)$ and $h(n)$
in $\Htilde$.  To compute $f(u)$ and $f(v)$, we simply wait for this edge
to appear in the enumeration of edges in $\Htilde$:  when it does, $\DD(\Gtilde)$
will have added a transcendental point $(\utilde,\vtilde)$ to the abelian variety
defined by $y^2=x(x-1)(x-f(a_m)-f(a_n))$ in $\Ktilde$, and we set
$f(u)=\utilde$ and $f(v)=\vtilde)$.  Then we extend $f$ to everything generated
by this $u$ and $v$ over the previous elements of $K_H$.  Finally,
when taking the differential closure $K$ itself of the field $K_H$ generated
by all of these pairs $(u,v)$ (in many distinct abelian varieties),
we need to be able to compute an isomorphism between the differential closures
(as given by Harrington) of isomorphic differential fields.
We now give the details of this last step, along with the earlier step
where details remain to be stated.

These two steps that require significant attention are analogous:  first, the extension of $f$
from $A$ to the algebraic closure $K_0$ of $\QA$;
and second, the final extension of $f$ from $K_H$ to its differential closure $K$.
The first of these is readily handled.  Since $A$ itself forms a
field-theoretic transcendence basis for the field $\QA$, we have
a splitting algorithm for $\QA$, and therefore can extend $f\res A$
effectively from this transcendence basis to the entire algebraic closure of $\QA$,
giving $f\res K_0$.

The effectiveness of the second step, extending the isomorphism
$f:K_H\to K_{\Htilde}$ between differential fields to an isomorphism
between their differential closures, is less well known.
Little work has been devoted to computable
categoricity for differential fields, as the uncertainty about whether $\TQbar$
is decidable has discouraged serious study.  Here we take some early steps
in that direction, using as an oracle the set $\TKHbar$ of constrained pairs
(equivalently, of complete formulas) over the differential field $K_H$.
(With this relativization we avoid the question of decidability of the set of constrained pairs!)
Corollary \ref{cor:oracle} below enables us to compute the oracle set
$\TKHbar$ uniformly from the oracle $T_\infty$ given for the reduction $\DD$.
We will not settle any major questions about computable categoricity here,
but we will do enough to compute the isomorphism $f=\DD(g)$
from $K$ onto $\Ktilde$ using that oracle.
%Notice that, since the differential field $K_H$ was produced effectively
%(uniformly in the oracle $H$), the set $\TKHbar$ is $\Pi^0_1$ uniformly
%in $H$, as remarked in Definition \ref{defn:constrainedpair}.

The tool essential for this process is a theorem developed in \cite{MOT14}
by Ovchinnikov, Trushin, and the author.

\begin{thm}[Theorem 8.6 of \cite{MOT14}]
\label{thm:algext}
Let $L$ be any computable differential field
with one nontrivial derivation, and $K$ its image under
any differential Rabin embedding of $L$ into any $\Khat$.
Then for every $z\in\Khat$, the constraint set $\overline{T_{K\la z\ra}}$
is computable in an oracle for $T_K$,
%Also, the constrainability set $\overline{U_{K\la z\ra}}$ is computably enumerable relative
%to $T_K$.  Both the computation of $T_{K\la z\ra}$ and the enumeration of $\overline{U_{K\la z\ra}}$
uniformly in $z$ and $T_K$.
\end{thm}

(Those authors also proved, as Theorem 9.6 of \cite{MOT14}, that the same holds
for $\overline{T_{K\la z\ra}}$ when $z$ is differentially transcendental over $K$.
They left open the analogous question for $z$ realizing any other nonprincipal type
over $K$, and to our knowledge that question remains open as of this writing.)

\comment{
\begin{thm}[Theorem 9.6 of \cite{MOT14}]
\label{thm:transext}
Consider any computable nonconstant differential field $K$, and
assume that $\Khat$ is not an algebraic field extension of $K$.
Let $K\la z\ra$ be
a computable differential field extension generated by an element
$z$ which is differentially transcendental over $K$, presented so that
$K$ is a computably enumerable subset of $K\la z\ra$.
Then the constraint set $\overline{T_{K\la z\ra}}$
is computable in an oracle for $T_K$,
%Also, the constrainability set $\overline{U_{K\la z\ra}}$ is computably enumerable relative
%to $T_K$.  Both the computation of $T_{K\la z\ra}$
%and the enumeration of $\overline{U_{K\la z\ra}}$
uniformly in $z$ and $T_K$.
\end{thm}
}

Actually, we require a modest generalization
of this theorem:  it must hold
not only for computable differential fields $L$ and $K$, but also for arbitrary
countable differential fields, assuming that we are given the atomic diagrams
$D(L)$ and $D(K)$ as oracles.  Moreover, the uniformity also carries over:
the computation of $\overline{T_{K\la z\ra}}$ in Theorem \ref{thm:algext}
can be carried out by a single Turing functional, uniformly in $z$ and the oracles
$T_K$ and $D(K)$ (which in turn are uniformly computable from $T_L$ and $D(L)$).
%while in Theorem \ref{thm:transext}, another Turing functional computes
%$\overline{T_{K\la z\ra}}$ uniformly in $z$ and the oracles $T_K$ and $D(K)$.
An examination of the proof of the original theorem will reveal that this
version does hold, relativized and uniformized, using exactly the same procedure as in the original.

Our extension of $f\res K_H$ to its differential closure $K=\DD(G)$
begins by using Corollary \ref{cor:oracle} to compute $T_{K_H}$,
using its enumeration $H=\HH(G)$ of the graph in question.
In particular, it finds a constrained pair $(p,q)\in \TKHbar$
such that the very first element $y_0$ (in the domain $\omega$ of $K$,
denoted here by $\{ y_0,y_1,\ldots\}$) satisfies $p(y_0)=0\neq q(y_0)$.
Such a pair must exist, and when it is found, we define $f(y_0)$
to equal the least element $\ytilde_i$ in the domain $\{\ytilde_0,\ytilde_1,\ldots\}$
of $\Ktilde=\DD(\Gtilde)$ such that $\ptilde(\ytilde_i)=0\neq \qtilde(\ytilde_i)$.
(Notice that, with $f\res K_H$ already defined, we may move effectively
between a differential polynomial such as $p$ in $K_H\{ Y\}$ and
its image in $K_{\Htilde}\{ Y\}$, denoted $\ptilde$, by mapping the coefficients of $p$
to their $f$-images in $K_{\Htilde}$.)  Since $K_H\cong K_{\Htilde}$,
some such $\ytilde_i$ must exist, and the first one we find becomes
$f(y_0)$.

Now we continue extending $f$ to map $K$ onto $\Ktilde$ by going back and forth
between them.  At each step, $f$ is already defined on $K_H\cup\{ y_{i_1},\ldots,y_{i_m}\}$
for finitely many additional elements of $K$, and we apply Theorem \ref{thm:algext}
$m$ times to compute $T_{K_H\la y_{i_1},\ldots,y_{i_m}\ra}$ from $T_{K_H}$;
likewise for $\Ktilde$.  Continuing this back-and-forth construction,
we clearly succeed in building the desired isomorphism $f=\DD(g)$
from $K=\DD(G)$ onto $\Ktilde=\DD(\Gtilde)$, and our construction
is effective uniformly in $g$, $D(G)$, and $D(\Gtilde)$.

This completes our computation of the reduction $\DD$ and the functor $\GG$,
using the oracle set $T_\infty$ to compute $\DD$, but no oracle for $\GG$.
Next we wish to show that their compositions $\DD\circ\GG$ and $\GG\circ\DD$
are essentially the identity maps on their respective domains.  First, for
$\GG\circ\DD$, we actually have the identity.
\begin{prop}
\label{prop:GcircD}
For every countable graph $G$, $\GG(\DD(G))=G$.
(That is, these are not just isomorphic graphs; they are the exact same presentation.)
\end{prop}
\begin{pf}
The key here is simply that the functors consider all elements in their proper
order as numbers in the domain $\omega$.  The graph $H=\HH(G)$
uses coding nodes $c_0<c_1<c_2<\cdots$ to code the nodes $0,1,2,\ldots$ of $G$,
and then fills in other nodes around them to create $H$.  Then the differential field
$K=\DD(G)$ makes the elements $a_0<a_1<a_2<\cdots$ of the quantifier-free-definable
set $A$ correspond to the nodes $0,1,2,\ldots$ of $H$, which forces
$a_{c_0}<a_{c_1}<a_{c_2}<\cdots$.  In the reverse direction, $\GG$ can decide
membership in $A$ from the atomic diagram $D(K)$ and it takes care
to number the nodes of the graph it produces in the same order that they appear in $A$:
thus we recover the same $H$ as before.

The final step is the one that requires a moment's thought.  The nodes of $H$
are partitioned into five classes: $\{c_n\}$, $\{ d_{mn}\}$, and those contained
in loops of length $3$, $5$, or $7$.  These classes do not overlap, and each is
computably enumerable.  So, for each node $p\in\omega=\dom{H}$ that
is enumerated into the first class (because we found a loop of length $3$ attached to $p$),
the functor $\FF$ ensures that
the corresponding node in $\FF(H)$ is numbered $m$, where $m=|\set{c_n}{c_n < p}|$.
(We did not give this much detail about $\FF$ earlier, but now it is relevant.)
Thus, even though we have only an enumeration of the graph $H$, we still make
the node $m\in\FF(H)$ correspond to the correct coding node $c_m$ in $H$.
From this is it now clear that $\GG(\DD(G))$ is just $G$ itself.
\qed\end{pf}

In the opposite direction, we do not get the identity.  Now we must state
our proposition more carefully, and must use an oracle set $T_\infty$.
\begin{prop}
\label{prop:DcircG}
There exists a Turing functional $\Lambda$ such that,
for every object $K\in\DCFT$, the function $\Lambda^{T_\infty\oplus K}$
is an isomorphism from $K$ onto $\DD(\GG(K))$.
\end{prop}
For reductions of categories, such a functional is the analogue of an effective natural
isomorphism between functors (cf.\ \cite[Defn.\ 1.6]{HTM3}).  Here, however, it is only
$T_\infty$-computable, not effective (unless $T_\infty$ is decidable).
\begin{pf}
The same arguments as for Proposition \ref{prop:GcircD} seem to apply here,
and indeed they do show that, for every $a_n$ in the subset $A$ of the given $K$,
the ``image'' of $a_n$ under the codings of $\DD\circ\GG$ will indeed be $a_n$ itself.
However, the other elements of $K$ are impossible to preserve this way.
Instead, we start with the knowledge that our map
$\lambda=\Lambda^{T_\infty\oplus D(K)}$ must send each $a_n$ to itself,
and then extend $\lambda$ to a full isomorphism from $K$ onto $\DD(\GG(K))$
in exactly the same way that $\DD$ itself begins with a bijection $\DD(g)\res A$ from
$A$ (as a subset of $\DD(G)$) onto itself (as a subset of $\DD(\Gtilde)$,
where $g:G\to\Gtilde$) and extends it to an isomorphism from all of $\DD(G)$
onto $\DD(\Gtilde)$.  The process is identical, extending first to 
$K_0$, then to $K_H$, and finally to all of $K$, using Theorem \ref{thm:algext}.
\qed\end{pf}

\subsection{The Oracle Set $T_\infty$}
\label{subsec:oracle}

Here we give the technical details about the oracle set $T_\infty$
that has been used in the construction in this section.  For these purposes,
let $G_\infty$ be the complete graph on the domain $\omega$,
and $H_\infty$ an enumeration of its edges (namely, all pairs $(m,n)$ with $m\neq n$).
We write $K_\infty=K_{H_\infty}$ to denote the differential field built by the
preceding construction on $H_\infty$, whose differential closure is $\DD(G_\infty)$.
(Notice that officially $G_\infty\notin\Graphs$, as the complete graph is automorphically trivial.
However, it is still possible to run the foregoing construction $\DD$ on $G_\infty$.)
Now we can define $T_\infty = T_{K_\infty}$ to be the constraint set for the
differential field $K_\infty$.  This is the oracle set used above.  As remarked earlier,
this set $T_\infty$ is computably enumerable, but it is an open question
whether it is decidable, and we will not address that question here.

For us, the relevant property of $T_\infty$ is that it is universal among the constraint
sets $K_H$, as $H$ ranges over all graphs on the domain $\omega$, and that this
universality is uniform.  This is expressed formally in Corollary \ref{cor:oracle} below,
after we prove the two necessary propositions.

\begin{prop}
\label{prop:upwards}
Let $H$ be (an enumeration of the edges in) a graph on the domain $\omega$,
and $K_H$ the differential field built earlier from this $H$.
Suppose $(p,q)\in\overline{T_{K_{H}}}$.  
Let $F$ be (an enumeration of the edges in) a graph on $\omega$ such that
every edge of $H$ is an edge of $F$.  Define the differential field embedding
$d:K_H\embeds K_F$ with $a_n\mapsto a_n$ (so $f\res\QA$ is the identity)
and $d(u_{ij})=u_{ij}$ and $d(v_{ij})=v_{ij}$ for all edges $(i,j)\in H$.  Write
$p^d$ for the image of $p$ in $K_F\{ Z\}$ under the map $d$ on the
coefficients of $p$, and similarly for $q^d$.  Then $(p^d,q^d)$ lies in $\overline{T_{K_{F}}}$.
\end{prop}
The map $d$ need not be the identity map on the domain $\omega$, strictly speaking,
because an edge $(i,j)$ of $H$ may have been enumerated into $F$ at a different stage,
making $u_{ij}$ a different element of $\omega$.
Here, for simplicity, we will regard $K_H$ as a differential subfield of $K_F$ via $d$,
so that $p$ and $p^d$ denote the same differential polynomial, and likewise for $q$.
Thus the Proposition claims that $(p,q)$ lies in $\TKFbar$ whenever it lies in $\overline{T_{K_{H}}}$.
\begin{pf}
Fix $F$, $H$, and $p,q\in K_H\{ Z\}$ as described, and let $f\in K_H\{ Z\}$.
If $\DCF0\cup D(K_H)\proves \forall Z[p(Z)=0\neq q(Z) \to f(Z)=0]$, then
only some initial segment $\sigma$ of $D(K_H)$ is used in the deduction.
But there is an enumeration $E$ of the edges in $F$ which enumerates
all edges outside of $H$ so late that $D(K_E)$ begins with that same
initial segment $\sigma$.  Thus $\DCF0\cup D(K_E)\proves\forall Z[p(Z)=0\neq q(Z) \to f(Z)=0]$,
and since $K_E\cong K_F$, $f(Z)=0$ lies in the type generated by
$p=0\neq q$ over $K_F$ as well.

The same argument holds if $\DCF0\cup D(K_H)\proves\forall Z[p(Z)=0\neq q(Z) \to f(Z)\neq 0]$.
Thus $p=0\neq q$ determines the truth of $f=0$ for all $f\in K_H\{ Z\}$,
since $(p,q)\in\overline{T_{K_{H}}}$.  It remains to consider those $f\in K_F\{ Z\}$
whose coefficients do not all lie in $K_H$.  (In fact, the argument for these $f$ subsumes the cases above.)

We can view any single $f\in K_F\{ Z\}$ as a differential polynomial over $K_H$:
$$ f(Z) = g(\uvec,\vvec,Z),\text{~with~}g(\Uvec,\Vvec,Z) \in K_H\{ \Uvec,\Vvec,Z\}$$
for some finite collection $\set{(u_i,v_i)}{i\leq j}$ of transcendental points
of elliptic curves corresponding to edges $(m_i,n_i)$ in $F$ but not in $H$.
Expressing $f$ this way, with coefficients from $K_H$,
we claim that there cannot exist $z_0,z_1\in\KFhat$, both realizing $(p,q)$,
such that $g(\uvec,\vvec,z_0)=0\neq g(\uvec,\vvec,z_1)$.
Indeed, in $\KHhat$, the following holds:
\begin{align*}
(\forall Z_0,Z_1)~
(\forall\text{~distinct~}A_0,\ldots,A_n\in A)
(\forall U_0,V_0,\ldots,U_j,V_j)
\\
\left[ \begin{array}{c}
p(Z_0)=p(Z_1)=0\neq q(Z_0)q(Z_1) \implies \vspace{2mm}\\
\left[ \begin{array}{c}
[(\forall i\leq j)~V_i^2 = U_i(U_i-1)(U_i-A_{m_i}-A_{n_i})]\to\\
(g(\Uvec,\Vvec,Z_0)=0 \leftrightarrow g(\Uvec,\Vvec,Z_1)=0)\end{array}\right]
\end{array}\right]
\end{align*}
(This sentence uses the definition of the infinite set $A$ of indiscernibles.)
Notice that the truth of this statement in $\KHhat$ is not vacuous:  there do exist
pairs $(u_i',v_i')\in(\KHhat)^2$ satisfying the elliptic-curve equations,
but they are all algebraic over $\QA$, as $\KHhat$ contains no transcendental points
for these curves.  (The corresponding edges lie in $F$ but not in $H$.)
Consequently all such coordinates $u_i'$ and $v_i'$ in $\KHhat$ lie in $\QAbar$, hence in $K_H$,
so that $f(\vec{u'},\vec{v'},Z)\in K_H\{ Z\}$.  With $(p,q)\in\overline{T_{K_{H}}}$,
the truth of the statement in $\KHhat$ now follows.

Since this sentence
(and its quantifier-free equivalent under \DCF0) hold in $\widehat{K_H}$,
it must have been proven by $\DCF0\cup\sigma$,
for some finite initial segment $\sigma$ of $D(K_H)$.
Now we choose an isomorphic copy $J$ of $K_F$ such that $\sigma\subset D(J)$,
which is possible since $K_H$ embeds into $K_F$ (as all edges in $H$ lie in $F$).
%% COULD TRY TO DO THIS EVEN IF $H$ HAS EDGES NOT IN $F$:
%% BUILD $K_F$ SO THAT TRANSCENDENTAL POINTS FOR SUCH EDGES IN $H$
%% BECOME TORSION POINTS IN $K_F$.  THIS REQUIRES MESSING WITH $K_0$ FURTHER OUT.....
%% STILL, IF IT WORKS, IT SHOWS THAT ALL $K_H$ HAVE TURING-EQUIV $T_{K_H}$.
Then $\DCF0\cup D(J)$ also proves the sentence above, which must therefore hold
in the differential closure of $J$, hence also in $\KFhat$.
Therefore we may plug in the elements $u_i$ and $v_i$ of $K_F$
for each $U_i$ and $V_i$, and the corresponding $a_m$ for each $A_m$,
giving the conclusion in $\KFhat$:
$$(\forall Z_0,Z_1)~[ p(Z_0)=p(Z_1)=0\neq q(Z_0)q(Z_1) \implies (f(Z_0)=0 \leftrightarrow f(Z_1)=0)].$$
With this, we see that even for those $f\in K_F\{ Z\}$ outside $K_H\{ Z\}$,
the formula $p(Z)=0\neq q(Z)$ does determine whether $f(Z)=0$ in $\KFhat$,
and so indeed $(p,q)\in\TKFbar$.
\qed\end{pf}

We also have a converse of Proposition \ref{prop:upwards}.
\begin{prop}
\label{prop:downwards}
As in Proposition \ref{prop:upwards}, let $F$ and $H$ be enumerations
of the edges of two graphs on the domain $\omega$, such that every edge
in $H$ is also in $F$.  Fix any $p,q\in K_H\{ Y\}$.  If $(p,q)\in\TKFbar$,
then also $(p^d, q^d)\in\overline{T_{K_H}}$.
\end{prop}
\begin{pf}
Again we view $K_H$ as a subfield of $K_F$ via the embedding $d$.
There is a quick but dangerous argument for this proposition:
every $z_0$ and $z_1$ in $\KFhat$ realizing $p(Z)=0\neq q(Z)$ are zeroes
of exactly the same polynomials in $K_F\{ Y\}$, and therefore of exactly
the same polynomials in $K_H\{ Y\}$.  This much is true, but in order to infer
Proposition \ref{prop:downwards}, one would need to know that every $z_0\in\widehat{K_H}$
lies in $\KFhat$.  For differential fields in general this does not hold
(see, e.g., \cite[\S 6]{M06}), so it would be necessary to prove it
for the specific situation of this $K_H$ and $K_F$.

However, there is a legitimate and even quicker proof:
$\TKFbar$ is a $\Pi_1$ subset of $(K_F^{<\omega})^2$
(albeit definable only by an $L_{\omega_1\omega}$ $\Pi_1$ formula),
and $K_H$ is a substructure of $K_F$, so $\TKFbar\cap (K_H^{<\omega})^2\subseteq\TKHbar$.
\qed\end{pf}

\comment{
Instead, we give the following proof.
Again we view $K_H$ as a subfield of $K_F$ via the embedding $d$.
Suppose $(p,q)\notin \overline{T_{K_H}}$.  Then there exists a polynomial
$h\in K_H\{ Y\}$ such that $\DCF0\cup D(H)$ proves
$$(\exists Z_0,Z_1)~[p(Z_0)=p(Z_1)=0\neq q(Z_0)q(Z_1)~\&~h(Z_0)=0\neq h(Z_1)].$$
Once again, a finite initial segment $\sigma\subset D(H)$ is all that is
needed for the proof, and this $\sigma$ can then be extended to the
atomic diagram of a copy $J\cong K_F$.  Thus $\DCF0\cup D(J)$
proves the same statement, so it must hold in $\widehat{J}$,
hence also in $\KFhat$, making $(p,q)\notin \TKFbar$.
}

\begin{cor}
\label{cor:oracle}
Let $T_\infty=T_{K_{H_\infty}}$ be the constraint set for the computable
differential field $K_\infty=K_{H_\infty}$ built from the complete graph $G_{\infty}$, containing
every edge $(m,n)\in\omega^2$ with $m\neq n$.
There is a Turing functional $\Theta$ such that, for every enumeration $H$
of any countable symmetric irreflexive graph on the domain $\omega$,
$$ T_{K_H} = \Theta^{H\oplus T_\infty}.$$
Thus $T_\infty$ may be seen as universal among the constraint sets for the
differential fields $K_H$.  (This covers not just graphs in $\Graphs$,
but also all automorphically trivial graphs, including $G_\infty$ itself.)
\end{cor}
\begin{pf}
This follows directly from Propositions \ref{prop:upwards} and
\ref{prop:downwards}.  In the details (as expressed in the statement of
Proposition \ref{prop:upwards}), $H_\infty$ plays the role of $F$, and we can compute
the embedding $f:K_H\embeds K_\infty$ using the enumeration $H$:
for each $x_{mn}$ or $y_{mn}$ in $K_H$, just wait until the edge $(m,n)$ appears in $H_\infty$,
at which point the corresponding $x_{mn}$ and $y_{mn}$ will be defined in $K_D$.
From our oracle we can decide whether $(f\circ p,f\circ q)\in T_\infty$,
and thus whether $(p,q)\in\TKH$.
\qed\end{pf}

In fact, the use of the oracle $T_\infty$ here is rather sparing.  What we really have is an
$m$-reduction from every $\TKH$ to $T_\infty$, computable uniformly in $H$.
That is, we have a functional $\Lambda$ such that, for each $H$, $\Lambda^H$
is a total function, with domain $(K_H\{ Z\})^2$,
for which $(p,q)\in\TKH$ just if $\Lambda^H(p,q)\in T_\infty$.
If desired, this function can readily be made injective for each single $H$,
giving a uniform $1$-reduction.

\section{The algebraic dependence relation}
\label{sec:binarydependence}

The predicate $C$ chosen for use in the category $\DCFC$ may seem unnatural,
as algebraicity over the differential subfield $\QA$ is not often considered
in differential algebra.  It makes our constructions in this article more transparent,
and we will continue to use it in subsequent sections.  Here, however,
we show that it could be replaced by the binary algebraic-dependence
predicate $B(x,y)$, which holds of two
elements $x$ and $y$ in a differential field $K$ if and only if the subfield
$\Q(x,y)$ of $K$ has transcendence degree $<2$ over $\Q$.  (To be clear:
$\Q(x,y)$ is a subfield but not necessarily a differential subfield, and we are
discussing the transcendence degree in the language of fields, rather than
the differential transcendence degree.)  So $B$, like $C$, is defined by a computable
$\Sigma^0_1$ formula of $L_{\omega_1\omega}$.
The main result is as follows.

\begin{prop}
\label{prop:binarydependence}
There exist Turing functionals $\Phi$ and $\Psi$ such that:
\begin{itemize}
\item
For every $K\in\DCFC$ isomorphic to a structure of the form $\DD(G)$,
$\Phi^{T_\infty\oplus K}$ computes the binary dependence relation $B^K$ on $K$; and
\item
For every $K\in\DCFields$ isomorphic to (the reduct of) a structure of the form $\DD(G)$,
$\Psi^{T_\infty\oplus K\oplus B^K}$ computes the unary relation $C^K$ on $K$.
\end{itemize}
\end{prop}
\begin{pf}
First we describe $\Psi$.  On an input $z\in K$, $\Psi$ searches for
$(a_{i_1},\ldots,a_{i_{2n}})\in A^{<\omega}$, pairs $(x_{1},y_{1}),\ldots,(x_{n},y_{n})$
from $K$, and polynomials $p,q\in \QA\{X_1,Y_1,\ldots,X_n,Y_n,Z\}$ such that:
\begin{itemize}
\item
$(\forall j\leq n)~y_{j}^2 = x_{j}(x_{j}-1)(x_{j}-a_{i_{2j-1}}-a_{i_{2j}})$; and
\item
For each $j\leq n$, the pair $(x_{j}, a_{i_{2j-1}}+a_{i_{2j}})$ does not lie in $B^K$; and
\item
$( p(x_{1},y_{1},\ldots,x_{n},y_{n},Z), q(x_{1},y_{1},\ldots,x_{n},y_{n},Z))\notin T_{K_H}$; and
\item
$p(x_{1},y_{1},\ldots,x_{n},y_{n},z) = 0 \neq q(x_{1},y_{1},\ldots,x_{n},y_{n},z).$
\end{itemize}
This $z$ must realize a principal type over the differential subfield $K_H$
(generated by the algebraic closure of $\QA$ along with all transcendental solutions
to the elliptic-curve equations over sums $(a_l+a_m)$ of elements of $A$),
as it lies in the differential closure $K$ of that $K_H$.  The procedure above simply finds
a constrained pair generating the type of $z$ over that subfield, using $B^K$ to recognize
transcendental solutions to the elliptic curves and using $T_\infty$ to compute the set $T_{K_H}$
as in Corollary \ref{cor:oracle}.  Edges in $H$ are enumerated as the list
of pairs $(l,m)$ for which a transcendental solution to the elliptic curve
on $(a_l+a_m)$ appears in $K$.  Eventually this search will locate the complete
formula generating the type realized by $z$ in $K$ over this $K_H$.
If $n=0$ (so no transcendental solutions were needed) and the $p$
it finds has order $0$ as a differential polynomial in $Z$, then $z$ is algebraic
over $\QA$ (as $p(z)=0$ and $p$ is an algebraic polynomial), so $C^K(z)$ holds in $K$.
If $p$ has order $0$ but $n>0$, then $z$ must be transcendental over $\QA$,
as $\{ x_{n},\ldots,x_{n}\}$ is algebraically independent over $\QA$ (and each
$y_{j}$ is interalgebraic with $x_{j}$), and so $C^K(z)$ fails in $K$.
Finally, if $p$ has positive order, then $z$ must be independent over
$\QA(x_{1},\ldots,y_{n})$, and so clearly $C^K(z)$ fails in $K$.

Next we turn to $\Phi$, which has access to $T_\infty$ and the atomic diagram of $K$
as a structure in $\DCFC$, meaning that it knows $C^K$.  Thus it can enumerate
transcendental solutions to elliptic-curve equations over sums of pairs of elements in $A$,
(much as $\Psi$ could, but now using $C^K$ instead of $B^K$), and thus can enumerate
both the graph $H$ and the corresponding differential subfield $K_H$
of which $K$ is a differential closure.  Given arbitrary $x\neq y$ in $K$, 
$\Phi$ first wishes to find the constrained pair $(p(X),q(X))\in\TKHbar$ realized by $x$ in $K$.
It does so using Corollary \ref{cor:oracle}.  If $p$ lies in $\Q[X]$, then immediately
it knows that $(x,y)\in B^K$.  (It is important here that $\Q$ be decidable within $K_H$,
but this is the content of Lemma \ref{lemma:splittingalg}.)

If $p\notin\Q[X]$, then $x$ is transcendental over $\Q$, and we continue by finding
a constrained pair over $K_H\la x\ra$ satisfied by $y$:  this process is the same as for $x$
except that it also requires the use of Theorem \ref{thm:algext}.  Once again, if we find
that $y$ satisfies an algebraic polynomial over $K_H(x)$, then we conclude that
$(x,y)\in B^K$.  Now, however, if the polynomial we found has degree $>1$ or positive order,
we have reached the opposite conclusion:  that $(x,y)\notin B^K$.
\qed\end{pf}

The resulting corollary proves our claim that we could have used the predicate $B$
in place of $C$ throughout this construction and also throughout Section \ref{sec:consequences}.

\begin{cor}
\label{cor:binarydependence}
There are $T_\infty$-computable functors in both directions, which are inverses of each other,
between the category 
$$ \set{K\in\DCFC}{(\exists G\in\Graphs)~K\cong\DD(G)}$$
and the category containing all the same countable differentially closed fields (on the domain $\omega$)
in the signature of differential fields without $C$ but with the binary-dependence predicate $B$
adjoined to the language and defined as above.  (In each category, the morphisms
are the differential field isomorphisms between objects.)
\end{cor}
\begin{pf}
The two functors simply map each structure in one category to the same structure
in the other category, and map each isomorphism to itself.  To be computable,
the functors merely need to compute $B^K$ from $C^K$ and vice versa,
which is the content of Proposition \ref{prop:binarydependence}.
\qed\end{pf}

\section{Consequences}
\label{sec:consequences}

Here we describe the consequences of the constructions in Section \ref{sec:functors}
for computable structure theory.  We will need to work ``on a cone,'' i.e.,
relativizing everything to the set $T_\infty$.

\subsection{Spectra}
\label{subsec:spec}

The first result is a refinement
of the main theorem from \cite{MM18}, quoted above as Theorem \ref{thm:MM}.
\begin{thm}
\label{thm:spectra}
On the cone above $T_\infty$,
the spectra of structures in $\DCFC$ are precisely the spectra of graphs
in $\Graphs$. For every $G\in\Graphs$, $\DD(G)$ satisfies
$ \text{Spec}_{T_\infty}(G) = \text{Spec}_{T_\infty}(\DD(G))$, i.e.,
$$  \set{ \deg{T_\infty\oplus D(\Gtilde)}}{\Gtilde\cong G} =
\set{ \deg{T_\infty\oplus D(\Ktilde)}}{\Ktilde\cong \DD(G)}.$$
Thus, on the cone above $T_\infty$, every spectrum of an automorphically nontrivial
structure in a computable language is realized as the spectrum
of a structure in $\DCFC$.
\end{thm}
\begin{pf}
It is well known that every spectrum of a structure in $\DCFC$ can be realized
as the spectrum of a graph; this is explained in \cite{HKSS02}.  (Clearly
all models of $\DCF0$ are automorphically nontrivial, and the
$L_{\omega_1\omega}$-definable predicate $C$ does not change this.)

For the converse, we claim that for every $G\in\Graphs$, the structure
$K=\DD(G)$ has the same spectrum as $G$ on the cone above $T_\infty$.
Indeed, for every $\Ktilde\cong K$, $\GG(\Ktilde)$ is a graph computable from
$\Ktilde$ and isomorphic to $\GG(K)=\GG(\DD(G))$,
which in turn is isomorphic to $G$ by Proposition \ref{prop:GcircD}.
Conversely, for every $\Gtilde\cong G$, $\DD(\Gtilde)$ is a copy of $\DD(G)$,
although now we can only claim that $\DD(\Gtilde)$ is computable from
$T_\infty\oplus\Gtilde$, since $\DD$ is $T_\infty$-computable.

Since $G$ and $K$ are both automorphically nontrivial, their $T_\infty$-spectra
are both upwards-closed under Turing reducibility (as seen in \cite{K86},
relativizing here to $T_\infty$),
and so these spectra must be equal.
\qed\end{pf}

\subsection{Computable categoricity}
\label{subsec:cc}

The remaining properties that we consider all involve isomorphisms
between graphs in $\Graphs$ or between structures in $\DCFT$.  Recall that,
although $\GG$ is a computable functor, $\DD$ is only a $T_\infty$-computable
reduction of categories, and it remains open whether the set $T_\infty$
is decidable or not.  Therefore, we usually work on the cone above $T_\infty$,
relativizing all statements to hold in a world where a decision
procedure for $T_\infty$ is given.  The first lemma is the key to essentially
all of our categoricity results.

\begin{lemma}
\label{lemma:iso}
Let $G,\Gtilde \in\Graphs$, and let $\bfd\geq \deg{T_\infty}$.
Then there is a $\bfd$-computable isomorphism between $G$ and $\Gtilde$
if and only if there is a $\bfd$-computable isomorphism between $\DD(G)$ and $\DD(\Gtilde)$.

Conversely, let $K,\Ktilde\in\DCFC$ be isomorphic to $\DD(G)$ (for the $G$ above),
still with $\bfd\geq \deg{T_\infty}$.
Then there is a $\bfd$-computable isomorphism between $K$ and $\Ktilde$
if and only if there is a $\bfd$-computable isomorphism between $\GG(K)$ and $\GG(\Ktilde)$.
\end{lemma}
\begin{pf}
If an isomorphism $g:G\to\Gtilde$ is $\bfd$-computable, then $\DD(g)$
is a $\bfd$-computable isomorphism from $\DD(G)$ onto $\DD(\Gtilde)$,
since $\bfd\geq\deg{T_\infty}$.  For the converse, for every
$\bfd$-computable isomorphism $f:\DD(G)\to\DD(\Gtilde)$, $\GG(f)$ is computable from $f$
(hence from $\bfd$) and, by Proposition \ref{prop:GcircD}, maps $G$ isomorphically
onto $\Gtilde$.

For the second part, $K$ and $\Ktilde$ both lie in the domain of the reduction
$\GG$, by assumption.  Therefore, an isomorphism $f:K\to\Ktilde$ will compute
the isomorphism $\GG(f):\GG(K)\to\GG(\Ktilde)$.  Conversely, if $g:\GG(K)\to\GG(\Ktilde)$
is a $\bfd$-computable isomorphism, then
$$ \left(\Lambda^{T_\infty\oplus \Ktilde}\right)^{-1} \circ \DD(g) \circ \Lambda^{T_\infty\oplus K}$$
is a $\bfd$-computable isomorphism $K\to\Ktilde$, where $\Lambda$ is as in 
Proposition \ref{prop:DcircG}.
\qed\end{pf}

For example, a $T_\infty$-computable
structure $\SS$ is \emph{computably categorical on the cone above $T_\infty$}
if and only if, for every $T_\infty$-computable $\widetilde\SS\cong\SS$, there exists
a $T_\infty$-computable isomorphism from $\SS$ onto $\widetilde\SS$.
The next definition generalizes this terminology.

\begin{defn}
\label{defn:compcatcone}
%Fix an oracle set $T\subseteq\omega$.  
Let $\A$ be a $T_\infty$-computable structure,
and let $\bfd\geq\deg{T_\infty}$ be a Turing degree.
Then $\A$ is \emph{$\bfd$-computably categorical on the cone above $T_\infty$}
if, for every $T_\infty$-computable $\B\cong \A$,
there is a $\bfd$-computable isomorphism from $\A$ onto $\B$.
\end{defn}

From Lemma \ref{lemma:iso}, we see first that on the cone above $T_\infty$,
$\bfd$-computable categoricity is preserved by $\DD$ and $\GG$.

\begin{cor}
\label{cor:cc}
Fix any degree $\bfd\geq\deg{T_\infty}$.
Then a $T_\infty$-computable graph $G$ is $\bfd$-computably categorical on the cone above $T_\infty$
if and only if $\DD(G)$ is.  Likewise, a $T_\infty$-computable $K\in\DCFC$
in the image of $\DD$ is $\bfd$-computably categorical on the cone above $T_\infty$
if and only if $\GG(K)$ is.  

Indeed, this statement remains true when $T_\infty$ is replaced everywhere in it
by any other oracle set $T\geq_T T_\infty$.
\qed\end{cor}

\begin{cor}
\label{cor:Pi11}
For $T_\infty$-computable structures in $\DCFT$, the property of being computably
categorical on the cone above $T_\infty$ is $\Pi^1_1$-complete.
\end{cor}
In contrast, for computable models of $\ACF0$, computable categoricity is only
$\Sigma^0_3$-hard:  such a structure is computably categorical if and only if it has
finite transcendence degree over $\Q$.  (The same holds for $T_\infty$-computable
categoricity for these structures:  it is $\Sigma^{T_\infty}_3$-hard.)
Moreover, for $\ACF0$, no additional definable predicates such as $C$ are necessary.
\begin{pf}
Since $T_\infty$ is $\Sigma^0_1$, the class $\Pi^1_1$ relative to $T_\infty$
is just the class $\Pi^1_1$, and the definition of $T_\infty$-computable categoricity
for $T_\infty$-computable structures in $\DCFC$ is clearly $\Pi^1_1$.  That it is
$\Pi^1_1$-hard follows from
Theorem 1 of \cite{DKLLMT15}, proven by Downey, Kach, Lempp, Lewis-Pye, Montalb\'an,
and Turetsky, which established that computable categoricity is a $\Pi^1_1$-complete property
of computable graphs.  Relativizing their argument to a fixed oracle such as $T_\infty$ is straightforward.
By Corollary \ref{cor:cc}, $T_\infty$-computable categoricity for $T_\infty$-computable graphs
is $1$-reducible to the same question for $\DCFC$.  To see this, notice that there is a
\emph{computable} injective function $h$ -- not merely $T_\infty$-computable --
such that, whenever $\Phi_e^{T_\infty}$ computes the atomic diagram of a graph $G$,
$\Phi_{h(e)}^{T_\infty}$ computes the atomic diagram of $\DD(G)$.  In particular,
$\Phi_{h(e)}$ uses its oracle $T_\infty$ both to compute the atomic diagram of $G$
(uniformly in the index $e$) and to compute the output of the $T_\infty$-computable
functor $\DD$ on that atomic diagram.
\qed\end{pf}

\subsection{Categoricity spectra}
\label{subsec:catspec}
\hyphenation{cate-gor-icity}

Lemma \ref{lemma:iso} also establishes that $\DD$ and $\GG$ preserve categoricity
spectra on the cone above $T_\infty$.  The relevant definition first appeared in \cite{FKM08};
its relativization to the cone above an arbitrary oracle set $T$ is used in \cite{CHT17}.
\begin{defn}
\label{defn:relcatspec}
For an oracle set $T$ and a $T$-computable structure $\A$,
the \emph{categoricity spectrum of $\A$ on the cone above $T$}
is the set of all Turing degrees $\bfd\geq\deg{T}$ such that $\A$ is $\bfd$-computably
categorical on the cone above $T$, as in Definition \ref{defn:compcatcone}.
If this set contains a least degree, then that degree is the \emph{degree of categoricity
of $\A$ on the cone above $T$}.
\end{defn}
\begin{prop}
For every $T_\infty$-computable structure $\SS$, there exists some
$K\in\DCFC$ with the same categoricity spectrum on the cone above $T_\infty$ as $\SS$.
\end{prop}
\begin{pf}
Hirschfeldt, Khoussainov, Shore, and Slinko proved in \cite{HKSS02} that for
automorphically nontrivial structures $\SS$, there is a graph $G\in\Graphs$ with the same categoricity
spectrum as $\SS$ on that cone, and by Corollary \ref{cor:cc}, the same holds of $\DD(G)$.

The categoricity spectrum of an automorphically trivial structure $\A$ always
contains all Turing degrees.  Taking $K=\DD(G_\infty)$ (for the complete graph $G_\infty$)
yields the same categoricity spectrum on the cone above $T_\infty$, although it is
unclear whether $\DD(G_\infty)$ is $\bfd$-computably categorical for degrees
$\bfd$ outside this cone.
\qed\end{pf}
If $T_\infty$ should turn out to be computable, this would also imply that computable
structures in $\DCFC$ realize all degrees of categoricity (as defined in \cite{FKM08})
as computable structures in general.  As it is, this clearly holds on the cone above $T_\infty$.

\subsection{Finite computable dimension}
\label{subsec:findim}

\begin{thm}
\label{thm:finitedim}
For every finite $n\geq 1$, there exists a $T_\infty$-computable structure in $\DCFT$
with exactly $n$ $T_\infty$-computable copies up to $T_\infty$-computable isomorphism.
\end{thm}
\begin{pf}
Fix any finite $n\geq 1$.
By results of Goncharov in \cite{G77, G80}, relativized to $T_\infty$,
there is a $T_\infty$-computable graph $G$
with exactly $n$ $T_\infty$-computable copies up to $T_\infty$-computable isomorphism.  Let
these be $G_1,\ldots,G_n$, and let $K_i=\DD(G_i)$, which is also $T_\infty$-computable.
By Lemma \ref{lemma:iso}, if $K_i$ and $K_j$ were $T_\infty$-computably isomorphic,
then $\GG(K_i)=G_i$ and $\GG(K_j)=G_j$ would be too, which is impossible for $i\neq j$.
So $K_1$ has at least $n$ copies, each $T_\infty$-computable, that are
pairwise not $T_\infty$-computably isomorphic.  

On the other hand,
if $\Ktilde$ is a $T_\infty$-computable copy of $K_1$, then $\GG(\Ktilde)$
is a $T_\infty$-computable copy of $G$, hence isomorphic to some $G_i$
via some $g\leq_T T_\infty$.  Now $\DD(g)$ will be a $T_\infty$-computable
isomorphism from $\DD(\GG(\Ktilde))$ onto $\DD(G_i)=K_i$.  By Proposition
\ref{prop:DcircG}, $\DD(\GG(\Ktilde))$ is $T_\infty$-computably isomorphic
to $\Ktilde$ itself, so $\Ktilde$ lies in the $T_\infty$-computable isomorphism class
of $K_i$.  This shows that $K_1$ has at most $n$ distinct $T_\infty$-computable copies
up to $T_\infty$-computable isomorphism, as required.
\qed\end{pf}

\subsection{Relative computable categoricity}
\label{subsec:rcc}

\begin{defn}
\label{defn:relrcc}
For an oracle set $S\subseteq\omega$, a countable structure $\A$ is
\emph{relatively $S$-computably categorical} if, for every pair
of structures $\B\cong\C\cong\A$, there exists an isomorphism
$f:\B\to\C$ which is computable from the oracle $S\oplus \B\oplus \C$.

More generally, for an $S$-computable ordinal $\alpha$, $\A$ is
\emph{relatively $\Delta^0_\alpha$-categorical above $S$} if, for every pair
of structures $\B\cong\C\cong\A$, there exists an isomorphism
$f:\B\to\C$ which is computable from the oracle $(S\oplus \B\oplus \C)^{(\alpha)}$,
the $\alpha$-th jump of the earlier oracle..
\end{defn}

Relative $S$-computable
categoricity is equivalent to the existence of a finite tuple $\avec$ from $\A$ and
a Scott family of existential formulas for $(\A,\avec)$
that has an $S$-computably-enumerable $e$-reduction to the existential
theory $\text{Th}_{\exists}(\A,\avec)$.
This generalizes the usual concept:  if $S$ itself is computable,
then such an $\A$ is said to \emph{relatively computably categorical} and the Scott
family is $e$-reducible to $\text{Th}_{\exists}(\A,\avec)$;
and if in addition $\A$ can be chosen to be computable,
then $\text{Th}_{\exists}(\A,\avec)$ is computably enumerable,
and thus so is the Scott family.  We refer the reader to \cite{M17}
for precise details.  Ash, Knight, Manasse, and Slaman originated
the concept in \cite{AKMS89}, and \cite{Mvol1,Mvol2} are useful sources as well.

\begin{thm}
\label{thm:rcc}
A graph $G$ is relatively $T_\infty$-computably categorical
if and only if $\DD(G)$ is.
\end{thm}
\begin{pf}
Let $\DD(G)$ be relatively $T_\infty$-computably categorical,
and suppose $G_0\cong G_1\cong G$.  Then there is some $\DD(G_0)\oplus\DD(G_1)$-computable
isomorphism $f:\DD(G_0)\to\DD(G_1)$, and $\GG(f)$ will map
$G_0=\GG(\DD(G_0))$ onto $G_1=\GG(\DD(G_1)) $ and will be computable
from $f$, hence from $\DD(G_0)\oplus\DD(G_1)$, hence from
$T_\infty\oplus G_0\oplus G_1$.  Thus $G$ itself is also
relatively $T_\infty$-computably categorical.  For the reverse implication
(when $G$ is relatively $T_\infty$-computably categorical),
one also needs Proposition \ref{prop:DcircG}:  if $K_0\cong K_1\cong \DD(G)$,
then each $g:\GG(K_0)\to\GG(K_1)$ yields an isomorphism
$$ \left(\Lambda^{T_\infty\oplus K_1}\right)^{-1} \circ\DD(g) \circ\Lambda^{T_\infty\oplus K_0}$$
mapping $K_0$ onto $K_1$ using $g$ and the oracle shown.
\qed\end{pf}

\begin{cor}
\label{cor:rcc}
Let $S\subseteq\omega$ satisfy $T_\infty\leq_T S$, and let $\alpha$ be an $S$-computable ordinal.
Then a graph $G$ is relatively $\Delta^0_\alpha$-categorical above $S$
if and only if $\DD(G)$ is.
\end{cor}
\begin{pf}
This is not strictly a corollary:  rather, for the proof, one repeats
the proof of Theorem \ref{thm:rcc}, taking $\alpha$-th jumps
in appropriate places.
\qed\end{pf}

\subsection{Uniform computable categoricity}
\label{subsec:ucc}

\begin{defn}
\label{defn:relucc}
For an oracle set $S\subseteq\omega$, a countable structure $\A$ is
\emph{uniformly $S$-computably categorical} if there exists a Turing functional
$\Gamma$ such that, for every pair of structures $\B\cong\C\cong\A$,
the function $\Gamma^{S\oplus \B\oplus \C}$ is an isomorphism
from $\B$ onto $\C$.

More generally, for an $S$-computable ordinal $\alpha$, $\A$ is
\emph{uniformly $\Delta^0_\alpha$-categorical above $S$} if there exists a Turing functional
$\Gamma$ such that, for every pair of structures $\B\cong\C\cong\A$,
the function $\Gamma^{(S\oplus \B\oplus \C)^{(\alpha)}}$ is an isomorphism
from $\B$ onto $\C$.
\end{defn}

This concept is very closely tied to Definition \ref{defn:relrcc}.  In general across all
these flavors, relative $\Delta^0_\alpha$-categoricity of $\A$ is equivalent to the
existence of a finite tuple $\avec$ of elements of $\A$ such that $(\A,\avec)$
is uniformly $\Delta^0_\alpha$-categorical, and uniform $\Delta^0_{\alpha+1}$-categoricity
of $\A$ is equivalent to the existence of a Scott family of (infinitary) $\Sigma^0_\alpha$
formulas for $\A$ (except that one must be careful with ``$\alpha+1$,'' distinguishing
finite and infinite ordinals).  For the proposition here, the main point is simply that the
procedure used in Corollary \ref{cor:rcc} is effective.

\begin{prop}
\label{prop:ucc}
Fix any countable ordinal $\alpha$ and any $S\subseteq\omega$ that can compute both $\alpha$
and $T_\infty$.  Then for each graph $G\in\Graphs$, $G$ is uniformly $\Delta^0_\alpha$-categorical
above $S$ if and only if $\DD(G)$ is.
\end{prop}
\begin{pf}
Suppose first that $\DD(G)$ is uniformly $\Delta^0_\alpha$-categorical above $S$.
Given any two copies $G_0$ and $G_1$ of $G$, the procedure for computing
an isomorphism $G_0\to G_1$ from an $(S\oplus G_0\oplus G_1)^{(\alpha)}$-oracle
begins by producing $(\DD(G_0))^{(\alpha)}$ and
$(\DD(G_1))^{(\alpha)}$, both of which can be done uniformly from the oracles
$G_0^{(\alpha)}$ and $G_1^{(\alpha)}$ using $\DD$ and the oracle $T_\infty$
(which can be computed from $S$).  The uniform categoricity procedure for $\DD(G)$
is then applied -- again using $S$ -- to compute an isomorphism $f:\DD(G_0)\to\DD(G_1)$,
and $\GG(f)$ must then be an isomorphism from $\GG(\DD(G_0))$ onto $\GG(\DD(G_1))$.
By Proposition \ref{prop:GcircD}, this is the desired output.

The converse works exactly the same way with the roles of the graph and
the differential field reversed, except that Proposition \ref{prop:DcircG}
is now required, as neither $\DD(\GG(K_i))$ is identical to the copy $K_i$ of $\DD(G)$
that we started with.  To fix this, we simply use $\Lambda^{T_\infty\oplus K_0}$
to map $K_0$ isomorphically onto $\DD(\GG(K_0))$, imitate the procedure above
to produce an isomorphism from $\DD(\GG(K_0))$ onto $\DD(\GG(K_1))$, and
finish off with $(\Lambda^{T_\infty\oplus K_1})^{-1}$, which maps $\DD(\GG(K_1))$
isomorphically onto $K_1$.
\qed\end{pf}

\subsection{Scott rank and categoricity ordinals}
\label{subsec:catord}

We will not take the time here to go through the various definitions of Scott rank,
but refer the reader instead to \cite{M15, Mvol1}.

Scott ranks raise some interesting questions here, as they come in both
parameterized and parameter-free versions.  When $\SS$ is a Scott family
of a graph and $\SS$ uses a parameter $c$ naming a node in the graph,
it seems natural to expect that if $n=c^G$ in the domain
$\omega$ of a particular presentation $G$ of the graph, then the corresponding
Scott family for $\DD(G)$ should require the element $a_n$
of the definable set $A$ in $\DD(G)$ to be named as a constant,
since $a_n$ ``represents'' $n$ in the interpretation of $G$ in $\DD(G)$.
Conversely, if an element of the definable set $A$ were used as a constant $d$
in a Scott family for some $K\in\DCFC$ from the range of $\DD$, and if
$a_m=d^K$ in a particular presentation of $K$, then the node $m$
of the graph $\GG(K)$ should be named as a constant in the corresponding Scott family
for $\GG(K)$.  However, it is also plausible that an element of $K$ not in $A$
might be used as a constant in a Scott family for $K$, and it is not clear
what constant(s) one might have to name in $\GG(K)$ to account for this.

Rather than address these issues here, we avoid them by appealing to the following
results about Scott rank, which appear in \cite{Mvol1}.  Our first notion, also
sometimes called the degree of categoricity, originated in \cite{CHT17}.

\begin{defn}
\label{defn:catordinal}
The \emph{categoricity ordinal} of a countable structure $\A$ is the least ordinal $\alpha$
for which there exists some $S\subseteq\omega$ such that $\A$
is relatively $\Delta^0_\alpha$-categorical above $S$.
\end{defn}

The term ``categoricity ordinal'' is not widely used, because by
Corollary VII.24 of \cite{Mvol2}, it is precisely the parameterized Scott rank
of the structure.  Thus, using our previous work,
we get a quick proof of the preservation of Scott ranks.

\begin{prop}
\label{prop:pScottrank}
For each $G\in\Graphs$, $\DD(G)$ has the same parametrized Scott rank as $G$.
Consequently, the structures in $\DCFC$ realize all possible parametrized Scott ranks:
every countable ordinal is the parametrized Scott rank of some $K\in\DCFC$.
\end{prop}
\begin{pf}
This is clear from Corollary \ref{cor:rcc}.  The fact that graphs realize all possible
parametrized Scott ranks is well established, e.g.\ by the construction in \cite{HKSS02}.
\qed\end{pf}

\cite[Corollary VII.24]{Mvol2} also establishes that the parameterless Scott rank
of a countable structure $\A$ is the least ordinal $\alpha$ such that there exists
some $S\subseteq\omega$ (that can compute $\alpha$) for which $\A$ is
uniformly $\Delta^0_\alpha$-categorical above $S$.  Hence we may apply
Proposition \ref{prop:ucc} in exactly the same way as Corollary \ref{cor:rcc}
to show that $\DD$ preserves these Scott ranks as well.

\begin{prop}
\label{prop:plessScottrank}
For each $G\in\Graphs$, $\DD(G)$ has the same parameter-free Scott rank as $G$.
Consequently, the structures in $\DCFC$ realize all possible parameter-free Scott ranks:
every countable ordinal is the parameter-free Scott rank of some $K\in\DCFC$.
\qed\end{prop}

This approach is evidently far simpler than analyzing finite tuples of constants directly
in the manner that was described (but not advised!) earlier in this subsection.

\section{Non-functoriality}
\label{sec:notafunctor}

The reduction $\GG$ constructed in Section \ref{sec:functors} is a natural
example of a functor, whose domain is the image of the reduction $\DD$.
However, $\DD$ is not a functor, in spite of all the efforts put into its construction.
It lacks the basic property known as \emph{functoriality}:  $\DD$ does not
in general respect the composition of isomorphisms $g_0:G_0\to G_1$
and $g_1:G_1\to G_2$.  The map $\DD(g_1\circ g_0)$ will indeed be an
isomorphism from $\DD(G_0)$ onto $\DD(G_2)$, but it will often fail to
equal $\DD(g_1)\circ\DD(g_0)$.  The reason for this failure lies in the
problem of extending $\DD(g_0)$ from $K_H$ to its differential closure,
where it is repeatedly required that we choose an element realizing a particular
type over the preceding elements.  We already made a concession to this problem
by allowing ourselves access to a $T_\infty$-oracle, simply in order to know
the types themselves.  The remaining problem is that many of these types
are realized by several elements -- indeed, often by infinitely many distinct
elements -- over the image of $\DD(g_0)$ defined so far, and there is no natural way to
choose among them.  We resorted to the expedient of using the natural ordering
of the domain $\omega$ and choosing the least element in that ordering that
realizes the type.  However, as that ordering has nothing to do with the structure
of $\DD(G_1)$, this destroys functoriality.

\subsection{Automorphism groups}
\label{subsec:autgp}

In fact, the quest for a pair of Borel functors $\DDtilde$
and $\GGtilde$, inverse to each other (up to a Borel natural transformation),
that would map $\Graphs$ into $\DCFields$ and back was doomed from the start.
As often happens, this is difficult to see from any of the discussion up till now,
but suddenly becomes clear when one considers automorphisms of the structures
in question.  David Marker was the first to realize that one should simply
ignore the differential structure.  The following argument, which for fields is folklore by now,
arose in discussions with him, Matthew Harrison-Trainor, and Tom Scanlon.

\begin{lemma}
\label{lemma:Marker}
No differentially closed field of characteristic $0$ has any automorphism of finite order $>2$.
Hence $\DCFields$ is not universal for automorphism groups.
\end{lemma}
\begin{pf}
Suppose that $K\in\DCFields$ has an automorphism $\alpha$ of finite order $n$.
Then $\alpha$ is also an automorphism of the reduct $K_0$ of $K$
in the language of fields (deleting the differentiation
operator), and its order remains the same.  We wish to apply the theorem of Artin
and Schreier (\cite[Satz 4]{AS65}, or see \cite{J89})
that no algebraically closed field has finite degree $>2$
over any subfield.  Indeed, if $F_0\subset K_0$ is the fixed field of $\alpha$,
then for each $t\in K_0\setminus F_0$, we have $\alpha^m(t)=t$ for some
(least) $m$ dividing $n$.  Now the polynomial
$$ \prod_{0\leq i<m} (X-\alpha^i(t))$$
has as its coefficients the elementary symmetric polynomials in $\{ t,\alpha(t),\ldots,\alpha^{m-1}(t)\}$,
which all lie in $F_0$ since $\alpha$ maps this set onto itself.  Thus every $t\in K_0$ has
degree $\leq n$ over $F_0$, and so $K_0/F_0$ is an algebraic extension, whose degree is
$\leq n$ by the Primitive Element Theorem.  By the Artin-Schreier Theorem,
this degree must be $\leq 2$, and so $\alpha\circ\alpha$ is the identity.
\qed\end{pf}

It now follows that there cannot exist any Borel functors $\DDtilde$
and $\GGtilde$ as described above, inverse to each other up to a
Borel natural transformation.  
%Indeed, fix any single $G\in\Graphs$
%with an automorphism of finite order $>2$.
%Then for every $K\in\DCFields$, there can be no Borel adjoint equivalence
%of categories between $\Iso{G}$ and $\Iso{K}$.  (
Here we use the language of \cite{HTM2}.
The ``equivalence of categories'' consists of Borel functors that are inverse
to each other up to a Borel natural transformation, exactly as stated above.
The following result of Harrison-Trainor, Montalb\'an and the author shows that
if such an equivalence did exist, then the automorphism groups of $K$ and $G$
would be isomorphic, contradicting Lemma \ref{lemma:Marker}.
\begin{thm}[Theorem 3.4 from \cite{HTM2}]
\label{thm:HTM2}
Suppose the functors $\FF \colon \Iso{\B} \to \Iso{\A}$, $\eta$,
$\GG \colon \Iso{\A} \to \Iso{\B}$, and $\epsilon$
form a Borel adjoint equivalence of categories between $\Iso{\A}$ and $\Iso{\B}$ with $\FF(\B) = \A$.
Then $\FF$, restricted to $\Aut{\B}$, gives an isomorphism between $\Aut{\B}$ and $\Aut{\A}$.
\end{thm}

Thus not only do the desired uniform functors between the categories $\Graphs$
and $\DCFields$ fail to exist, but there are even individual graphs for which
no such functor exists.

\begin{cor}[Harrison-Trainor, Marker, Miller \& Scanlon]
\label{cor:nofunctor}
There is no Borel adjoint equivalence of categories between 
$\Graphs$ and $\DCFields$.
Indeed, for certain graphs $G$, there is no Borel adjoint equivalence of categories
between $\Iso{G}$ and any $\Iso{K}$ with $K\in\DCFields$.
\end{cor}
\begin{pf}
By Lemma \ref{lemma:Marker}, take any graph $G$ with an automorphism
of finite order $>2$.
\qed\end{pf}

Moreover, this also
means that there can be no Borel bi-interpretation
between the $G$ above and any $K\in\DCFields$, as \cite[Theorem 2.6]{HTM2}
also shows that such a bi-interpretation would yield isomorphic automorphism groups.
(The notion of \emph{interpretation} used there is a substantial
generalization of the usual model-theoretic notion of an interpretation
by finitary formulas:  in \cite{HTM2}, $L_{\omega_1\omega}$ formulas
may be used, and the domain of the interpretation may consist of
tuples of arbitrary finite length, rather than a fixed length.)  In particular,
while the Hrushovski-Sokolovi\'c construction gives interpretations in both directions
between $G$ and $\DD(G)$, neither one can be half of a bi-interpretation.

Before Marker arrived at the straightforward argument given above,
he and the author had been considering whether a model of \DCF0
could be rigid.  There certainly do exist rigid countable graphs
(i.e., graphs with no nontrivial automorphisms), whereas the possibility of a rigid differentially
closed field appeared very unlikely.  Models of \ACF0 certainly always have nontrivial
automorphisms.  However, after the author raised this question, Marker investigated
and proved that the unlikely situation actually can hold.  So the delay in
discovering the simple argument above had the positive effect of catalyzing
the following theorem.
\begin{thm}[Marker, Theorem 3.5 in \cite{Mta}]
\label{thm:Dave}
There exists a rigid countable differentially closed field of characteristic $0$.
\end{thm}
The field establishing this theorem has the property that it is not a differential closure
of any proper differential subfield of itself; indeed this property is implied by rigidity
for countable models of \DCF0.  In contrast, every algebraically closed field
is the algebraic closure of some proper subfield, indeed of a real-closed subfield.

Marker notes that his construction does not appear to generalize to produce 
non-rigid but almost-rigid models of \DCF0, i.e., non-rigid models having
finite automorphism groups.  He poses the questions of whether \DCF0 has
any countable non-rigid almost-rigid models.  (By his Lemma \ref{lemma:Marker},
the automorphism group would need to consist entirely of involutions, hence would be
abelian of the form $(\Z/(2))^k$.)

With help from Marker, however, we can give a different proof here that for our specific reduction
$\DD$ from Section \ref{sec:functors}, there is no $\GGtilde$ that can serve as
its inverse, even up to a natural transformation.  The proof here uses the following fact.
\begin{prop}[Marker, Proposition 1.1 in \cite{Mta}]
\label{prop:closure}
Let $K\in\DCFields$, and assume that $K$ is not differentially closed.  Then its
differential closure $\Khat$ is not rigid.
\end{prop}
Marker's field in Theorem \ref{thm:Dave} avoids this problem by not being the
differential closure of any proper differential subfield of itself.
\comment{
\begin{pf}
Fix a differential closure $\Khat$ of $K$.
Since $K$ is not differentially closed, its constraint set $\TKbar$ contains some pair
$(p,q)\in (K\{ Y\})^2$ for which $p$ has either positive order or else degree $>1$.
Thus the formula $p(Y)=0\neq q(Y)$ is realized by at least two elements $a\neq b$ of $\Khat$,
and there exists an isomorphism $f:K\la a\ra\to K\la b\ra$ of differential fields with
$f(a)=b$ and $f\res K=\text{id}_K$.  By Proposition 5 in \cite[\S 4]{K74},
$\Khat$ is a differential closure of $K\la a\ra$ and also of $K\la b\ra$.  Therefore,
Corollary 2.10 of \cite{M06} shows that $f$ must extend to a map from $\Khat$
(viewed as $\widehat{K\la a\ra}$) onto
itself (viewed now as $\widehat{K\la b\ra}$), and this extension is a nontrivial
automorphism of $\Khat$.
\qed\end{pf}
}
The specific reduction $\DD$ that we constructed in Section \ref{sec:functors}
built a differential field $K_H$ from an arbitrary enumeration of the graph $H=\HH(G)$,
and then applied Harrington's Theorem to construct the differential closure $\DD(G)$
of $K_H$.  Since $K_H$ itself was clearly not differentially closed, Proposition
\ref{prop:closure} shows that $\DD(G)$ is not rigid.  Consequently, for each rigid graph
$G$, $\aut{\DD(G)}$ is not isomorphic to $\aut{G}$, eliminating all hope of converting
the Hrushovski-Sokolovi\'c construction into inverse functors.

It may be noted that the map $G\mapsto K_{\HH(G)}$ is indeed a functor on $\Graphs$,
with isomorphisms $G_0\to G_1$ mapping to isomorphisms
$K_{\HH(G_0)}\to K_{\HH(G_1)}$ in the obvious way.  Moreover, it is computable
(even without the predicate $C$ or the oracle $T_\infty$), and has a computable
inverse functor from its image back to $\Graphs$.  Using this functor, one may carry
universality results about $\Graphs$ over to the larger category of countable differential fields
(not necessarily differentially closed!)
of characteristic $0$ in full.  However, it was already known that countable differential
fields form a universal category:  the category $\Fields_0$ is known to be universal
in all these respects, using the functors constructed in \cite{MPSS18}, and
every field becomes a differential field simply by giving it the trivial derivation $0$.
We point out that \cite{MPSS18} also established universality for $\Fields_p$
with $p>0$,
using earlier results of Fried and Koll\'ar in \cite{FK82}, and that universality
for the category of differential fields of characteristic $p$ follows from
the same trick of making a field into a differential field with derivative $0$ everywhere.

\comment{
Finally, we remark that by the results of Harrison-Trainor, Melnikov, Montalb\'an,
and the author in \cite{HTM3,HTM2}, questions about the existence of functors
-- whether computable or at any other level within the Borel hierarchy --
are equivalent to questions about the existence of
uniform interpretations between objects in the categories in question.
Here ``interpretations'' is construed using a broad meaning, under which
formulas of $L_{\omega_1\omega}$ may be used to define the interpretation
and elements of one structure may be represented by tuples of arbitrary
finite length in the other structure.  The questions here about the existence
of functors may therefore be rephrased as questions about such interpretations
between these classes of structures, and questions of bi-interpretability
may also be addressed thus.
}

\section{Further questions}
\label{sec:questions}

The preceding section included a natural question of Marker regarding
automorphism groups, but this is not the only question
arising out of the results here.  The most obvious question left open
in this article is the decidability of the oracle set $T_\infty$.
No matter what the answer to this question is, compelling results would follow.
If $T_\infty$ is decidable, then Section \ref{sec:consequences} would give
universality results about categoricity for $\DCFC$ in
all the usual senses of \cite{HKSS02} and others, without any relativization to a cone.
This would say in a strong way that the addition of the infinitarily-definable
predicate $C$ to the signature fills the gap in the computability properties
of models of $\DCF0$.  On the other hand, while the concept of working on a cone
of Turing degrees has seen substantial and interesting applications in recent years,
there is no known natural class of structures for which it has yet proven necessary.
If $T_\infty$ is undecidable, then the results here would make it the first example
of such a class.

Recall that \DCF0 is the theory of \emph{ordinary} differentially closed fields,
i.e., with just a single derivation in the signature.  Of course, the universality
results of Section \ref{sec:consequences} for the class $\DCFC$ imply
the same for the larger class of partial differentially closed
fields, with an arbitrary finite number of derivations (which are usually assumed
to commute with one another).  Nevertheless, it is well known that the general study of ordinary
differential equations becomes significantly more complex when one switches
to partial differential equations, and so one naturally asks whether this increased complexity
manifests itself in some computability-theoretic way when more derivations
are added to the signature.  If anything, the relevance of this question dims a bit
in light of the present article, which shows that even the ordinary case comes
close to being universal.  (Also, since a differentially closed field with several derivations
is algebraically closed, the arguments of Section \ref{sec:notafunctor} still preclude universality
for automorphism groups.)  Nevertheless, it seems possible that with additional derivations,
one might no longer require the predicate $C$ to achieve universality.  Alternatively,
the oracle set $T_\infty$ might have to be replaced by another oracle, possibly of
different Turing degree.  Any of these would emphasize the distinction between
ordinary and partial differential algebra, whereas if they proved false, then the distinction
would seem less dramatic.

%\parbox{4.7in}{
%{\sc
%\noindent
%Department of Mathematics, Statistics, \& Computer Science \hfill \\
%\hspace*{.1in} University of Illinois at Chicago M/C 249 \hfill \\
%\hspace*{.2in}  851 S.\ Morgan St. \hfill \\
%\hspace*{.3in}  Chicago, IL  60607 U.S.A. \hfill }\\
%\medskip
%\hspace*{.045in} {\it E-mail: }
%\texttt{marker\at {math.uic.edu}%\at {up.edu}
% }\hfill \\
%}\\
\comment{
\parbox{4.7in}{
{\sc
\noindent
Department of Mathematics \hfill \\
\hspace*{.1in}  Queens College -- C.U.N.Y. \hfill \\
\hspace*{.2in}  65-30 Kissena Blvd. \hfill \\
\hspace*{.3in}  Flushing, New York  11367 U.S.A. \hfill \\
Ph.D. Programs in Mathematics \& Computer Science \hfill \\
\hspace*{.1in}  C.U.N.Y.\ Graduate Center\hfill \\
\hspace*{.2in}  365 Fifth Avenue \hfill \\
\hspace*{.3in}  New York, New York  10016 U.S.A. \hfill}\\
\medskip
\hspace*{.045in} {\it E-mail: }
\texttt{Russell.Miller\at {qc.cuny.edu} }\hfill \\
}
}

\end{document}